\documentclass{custom_template}

\usepackage{amsmath}
\usepackage[english]{babel}
\usepackage[utf8]{inputenc}
\usepackage[T1]{fontenc}
\usepackage{lmodern}
\usepackage{graphicx}
\usepackage{combelow}
\usepackage{mathtools}
\usepackage{amsfonts, amssymb}
\usepackage{sidecap}
\usepackage{lipsum}

\usepackage[colorlinks,allcolors=black]{hyperref}

\usepackage{setspace}
\title{\vspace{-1cm}\centering\bfseries\semilarge Speed and stability of segregated waves in a pressure-based model of heterogeneous cell populations
}

\author{\centering \normalsize Carles Falc\'o$^{1,*}$\;, Rebecca M.  Crossley$^{1}$\;, Martina Conte$^{2}$\;, and Tommaso Lorenzi$^{2}$}

\date{}

\newcommand{\nhat}{{\hat{\boldsymbol{\mathrm{n}}}}}
\newcommand{\vn}{V_\mathrm{n}}
\newcommand{\er}{\hat{\boldsymbol{e}}_r}

\newcommand{\vv}{{\mathrm{v}}}
\newcommand{\uu}{{\mathrm{u}}}

\newcommand{\dd}{{\mathrm{d}}}
\newcommand{\ww}{{\mathrm{w}}}

\numberwithin{equation}{section}

\begin{document}

\maketitle\par\vspace{-1.3cm}
\begin{center}
{$^{1}$Mathematical Institute, University of Oxford, Oxford, United Kingdom
\\
$^{2}$Department of Mathematical Sciences "G.L. Lagrange", Politecnico di Torino, Turin, Italy}

{\vspace{.1cm}\small $^*$Correspondence: falcoigandia@maths.ox.ac.uk }

\vspace{.4cm}

\begin{minipage}{16cm}  \small We consider a minimal pressure-based model of heterogeneous cell populations consisting of proliferative and non-proliferative cells with different mobilities. The model is formulated as a system of reaction--cross--diffusion equations describing the spatio-temporal dynamics of the cell densities. The model is known to admit one-dimensional travelling wave solutions with strictly segregated components: non-proliferative cells occupy a finite region at the leading edge, while proliferative cells remain at the rear. However, the speed, parameter dependence, and stability of these waves remain poorly understood. In this work, we derive an almost explicit variational bound on the wave speed by reformulating the problem as a free-boundary problem for a generalised porous--Fisher equation. The estimates we obtain apply to general pressure laws and growth kinetics, agree closely with the results of numerical simulations, and become sharp in the incompressible limit, where we formally recover a fully explicit characterisation of the wave speed. We then analyse the stability of the waves to show that segregated waves are stable only when non-proliferative cells are more mobile than proliferative cells. Finally, motivated by numerical observations of finger-like protrusions, we investigate the stability of incompressible segregated circular waves through asymptotic shape-perturbation analysis. This yields explicit expressions for the pressure, interface velocity, and growth rates of angular modes, thereby making evident the destabilisation mechanisms that may lead to the emergence of fingering instability. Interestingly, we find that, in contrast with the one-dimensional case, the stability of such circular waves is not determined solely by the relative value of the mobility coefficients, and thus instabilities may arise irrespective of which cell type has the larger mobility.
\end{minipage}
\end{center}

\par\vspace{.2cm}


\section{Introduction}
Tissue development, wound healing, and tumour growth are all shaped by the complex interplay between pressure-dependent inhibition of cell proliferation, where cell division is arrested at sufficiently high densities~\cite{byrne2003modelling,byrne2009individual,drasdo2012modeling,ranft2010fluidization}, and mechanically regulated cell movement. The latter is often simplified and represented as cell movement down the gradient of cellular pressure, towards regions where cells are less compressed~\cite{byrne1997free,byrne2003modelling}. Consequently, pressure-based models are often used to capture the growth dynamics of various cell populations~\cite{roose2007mathematical}.

Building on the ideas presented in the seminal paper describing pressure-driven cell movement~\cite{greenspan1976growth} and subsequent extension papers~\cite{ambrosi2002closure,bresch2010computational,byrne2009individual,byrne2003modelling,byrne1997free,ciarletta2011radial,lowengrub2009nonlinear,preziosi2009multiphase,sherratt2001new}, pressure-based models are built on three main assumptions: cell movement can be described using an advective term, where the advective velocity is inversely proportional to the gradient of the cellular pressure; cell proliferation and death can be described using a reaction term, with pressure-dependent growth kinetics; and the cellular pressure can be defined as a function of the cell density through a density-based law, that is, the pressure law is formulated as a \emph{barotropic} relation. A prototypical example of these models is provided by the following nonlinear reaction--diffusion equation~\cite{byrne2009individual}:
\begin{subequations}\label{eq:p-bmod1pop}
    \begin{align}
        \partial_t\rho & = \mu\nabla\cdot\left(\rho\,\nabla p\right) + \rho \,G(p)\,,
\label{eq:p-bmod1pop_a}\\
          p & := \Pi(\rho) \,,\label{eq:p-bmod1pop_b}
    \end{align}
\end{subequations}
where $\rho(\mathbf{x},t)$ is the density of the cells at position $\mathbf{x}\in\mathbb{R}^n$ and time $t>0$, and $p(\mathbf{x},t)$ is the cellular pressure. The first term on the right-hand side of Eq.~\eqref{eq:p-bmod1pop_a} models cell movement and, in analogy with the Darcy’s law for fluid flow in porous media~\cite{darcy1856fontaines}, the parameter $\mu>0$ is referred to as the cell mobility coefficient. The mobility coefficient is inversely proportional to the permeability of the medium in which the cells are embedded (\textit{e.g.}, the extracellular matrix) and depends on the cells’ morphological and mechanical properties~\cite{ambrosi2002closure,byrne2009individual}. Moreover, the second term on the right-hand side of Eq.~\eqref{eq:p-bmod1pop_a} takes into account cell population growth kinetics (\textit{i.e.}, cell proliferation and death), and the function $G(p)$ is the net growth rate depending on the cell density through the pressure $p$ and is commonly assumed to be a non-increasing function that vanishes at a critical value, $\bar p>0$, known as the \emph{homeostatic pressure}. The homeostatic pressure accounts for the fact that cells will stop dividing if the pressure at their current position exceeds such a critical value~\cite{basan2009homeostatic,shraiman2005mechanical}. As such, $G(p)$ typically satisfies the following assumptions:
\begin{equation}
\label{ass:G}
G(0) < \infty \, , \quad G(\bar p) = 0 \, , \quad G' < 0 \, .
\end{equation}
Finally, the cellular pressure is defined through the density-based law $\Pi(\rho)$ in Eq.~\eqref{eq:p-bmod1pop_b}, which is commonly assumed to be non-decreasing~\cite{ambrosi2002closure,byrne2009individual}. A possible choice that has recently received increasing attention from the mathematical community---see for instance \cite{bubba2020hele,david2024incompressible,kim2016free,mellet2017hele,perthame2014hele} and references therein---is the power-law
\begin{equation}\label{eq: pl pressureori}
    \Pi(\rho) := K_{\gamma} \, \rho^{\gamma} \, , \quad \gamma \geq 1 \, ,
\end{equation}
which can be derived from microscopic models when accounting for volume exclusion. Different exponents appear depending on the cells' size, shape, and microscopic interaction rules~\cite{bruna2017diffusion, bakerAspectRatio2, falco2022random}. In Eq.~\eqref{eq: pl pressureori}, the parameter $\gamma$ provides a measure of the stiffness of the pressure law, and $K_{\gamma}>0$ is a scale factor such that $K_{\gamma} \to 1$ as $\gamma \to \infty$. As demonstrated in~\cite{perthame2014hele}, under the pressure law defined via Eq.~\eqref{eq: pl pressureori}, a link between models of the form of Eq.~\eqref{eq:p-bmod1pop} and models formulated as free-boundary problems, which have also been widely employed to study the growth of cell populations~\citep{friedman2015free}, can rigorously be established in the limit $\gamma \to \infty$. Such an asymptotic regime is usually referred to as the \emph{incompressible limit}, because it corresponds to mathematically approximating cells as an incompressible fluid within the framework of pressure-based models.

Models of the form of Eq.~\eqref{eq:p-bmod1pop} and its related variants have drawn interest from mathematicians and physicists alike. Firstly, for their ability to recapitulate key aspects of both tumour and tissue growth, but also for exhibiting travelling waves which display interesting features---see for instance~\citep{bertsch2015travelling,carrillo2024multipop,chaplain2020bridging,lorenzi2016interfaces,tang2014composite} and references therein.

\subsection{A minimal pressure-based model for the growth of heterogeneous cell populations} 
\label{sec: preliminaries}
An implicit assumption of the model provided by Eq.~\eqref{eq:p-bmod1pop} is phenotypic homogeneity, whereby cells in the population are assumed to be, as a first approximation, identical. However, heterogeneity is typically observed within cell populations, and has been found to play a pivotal role in tissue development, wound healing, and tumour growth~\cite{huang2009non,meacham2013tumour,rognoni2018skin,wang2022cellular}. 
A possible extension of the model in Eq.~\eqref{eq:p-bmod1pop} to populations comprising two types of cells, \textit{i.e.}, non-proliferative and proliferative cells with different mobility coefficients, is given by the following system of reaction--cross--diffusion equations~\cite{lorenzi2016interfaces}, which provides a minimal pressure-based model for the growth of heterogeneous cell populations:
\begin{subequations}\label{eq:full_model}
    \begin{align}
        \partial_t\rho & = \mu\nabla\cdot\left(\rho\,\nabla p\right) + \rho \,G(p)\,,
\\
        \partial_t\eta &  = \nu\nabla\cdot\left(\eta\,\nabla p\right)\,,
        \\
          p & := \Pi(\sigma) \,, \quad \sigma := \rho + \eta \,.
    \end{align}
\end{subequations}
Here, $\rho(\mathbf{x},t)$ and $\eta(\mathbf{x},t)$ are the densities of the proliferative and non-proliferative cells, respectively, at position $\mathbf{x}\in\mathbb{R}^n$ and time $t>0$, while the cellular pressure $p(\mathbf{x},t)$ is now defined as a function of the total cell density $\sigma(\mathbf{x},t):=\rho(\mathbf{x}, t)+\eta(\mathbf{x},t)$. In analogy with Eq.~\eqref{eq:p-bmod1pop}, the parameters $\mu>0$ and $\nu>0$ are the mobility coefficients of proliferative and non-proliferative cells, respectively, and the function $G(p)$ is the net growth rate of the density of proliferative cells, which satisfies the assumptions given by Eq.~\eqref{ass:G}.

\subsection{Previous results and open questions}
Despite its apparent simplicity, the model provided by Eqs.~\eqref{eq:full_model} exhibits rich behaviour. 
In one spatial dimension, under monotone density-based laws for the cellular pressure satisfying
\begin{equation}\label{ass:pressure}
    \Pi(0) = 0\, , \quad \Pi' > 0 \, , 
\end{equation}
Eqs.~\eqref{eq:full_model} admit travelling wave solutions with strictly segregated components, \textit{i.e.}, the densities of the two cell types have disjoint supports, whereby non-proliferative cells are confined to a finite region at the leading edge, where the pressure approaches zero, while proliferative cells are restricted to the rear of the wave, where the pressure saturates at the homeostatic value $\bar p$~\cite{chaplain2020bridging,lorenzi2016interfaces}. However, to the best of our knowledge, no explicit estimates on the wave speed have been found. Moreover, while the analytical construction of such travelling wave solutions does not depend on the relative value of the mobility coefficients $\mu$ and $\nu$, numerical simulations indicate that Eqs.~\eqref{eq:full_model} can support segregated travelling waves only when $\nu>\mu$, \textit{i.e.}, when the mobility coefficient of the non-proliferative cells is larger than that of the proliferative cells. In contrast, if $\nu<\mu$ then the non-proliferative cells are left behind by the proliferative cells and eventually form their own wave in isolation---\textit{i.e.}, after an initial transient, the density of the non-proliferative cells converges to a stationary profile, while the density of the proliferative cells exhibits travelling front-like behaviour. While these numerical observations support the idea that segregated travelling wave solutions wherein non-proliferative cells form the front of the wave followed by a bulk of proliferative cells are stable only if $\nu>\mu$, a mathematical formalisation of this idea has not yet been developed. 

In two spatial dimensions, when a density-based law for the cellular pressure of the type of Eq.~\eqref{eq: pl pressureori} is considered, \textit{i.e.}, 
\begin{equation}\label{eq: pl pressure}
    \Pi(\sigma) := K_{\gamma} \, \sigma^\gamma\,,\quad \gamma\geq 1\,, \quad K_{\gamma} := \dfrac{\gamma+1}{\gamma} \, , 
\end{equation}
numerical simulations indicate that if $\nu>\mu$ then Eqs.~\eqref{eq:full_model} can support circular waves wherein the leading edge of the wave consists entirely of non-proliferative cells~\cite{lorenzi2016interfaces}. In contrast, if $\nu < \mu$, then finger-like protrusions formed by proliferative cells that protrude through and displace the non-proliferative cells ahead emerge. The striking behaviour, which is reminiscent of Saffman--Taylor instabilities arising when a less viscous fluid displaces a more viscous one in a porous medium or Hele-Shaw cell~\cite{saffman1958penetration}, is particularly relevant in the context of tumour growth and morphogenesis. In fact, fingering instabilities are observed both at the interface between tumour cells and adipose tissue in breast cancer invasion~\cite{wang2012adipose} and at the interface between epithelial cells and the extracellular matrix in budding morphogenesis~\cite{wang2021budding}. 

The numerical evidence for the formation of such finger-like protrusions led Kim and Jong~\cite{kim2021interface} to consider Eqs.~\eqref{eq:full_model} in an almost radially symmetric setting in two spatial dimensions. Here, an inner interface separates proliferative cells from non-proliferative ones, while an outer interface separates the latter from the vacuum. By deriving evolution equations for the two interfaces, Kim and Jong proved well-posedness of the resulting problem when $\mu<\nu$ for nearly radial initial configurations. Interestingly, their analysis shows that, upon linearising the inner interface around a perfect circle, the angular modes of the interface satisfy a fractional diffusion equation in which the diffusivity depends on the sign of $\mu-\nu$. When $\mu>\nu$, the diffusivity is negative, suggesting instability that may result in the formation of finger-like protrusions.

\subsection{Summary of the main results in the paper} The results presented in this paper complete the picture for segregated travelling waves exhibited by Eqs.~\eqref{eq:full_model}, resolving the aforementioned outstanding questions.

On the one hand, we use a variational argument to derive almost fully explicit bounds on the speed of segregated travelling wave solutions. The variational method builds on the classical variational principle for reaction--diffusion equations developed by Benguria and Depassier~\cite{benguria1994validity, benguria1996speed, benguria1996variational, benguria2004minimal, stokes2024speed}, and on its recent extensions to moving boundary problems and reaction--diffusion systems~\cite{crossley2026optimalcontrolapproachnonlinear}. As part of this analysis, we derive wave speed estimates for a generalised version of the porous--Fisher model with a moving boundary, which was previously studied by Fadai and Simpson~\cite{fadai2020new} and arises from the model provided by Eqs.~\eqref{eq:full_model} and~\eqref{eq: pl pressure} in the travelling wave framework. We show that the obtained results yield bounds on the wave speed that are in excellent agreement with simulations and closely match the numerically estimated values. Furthermore, we provide an explicit characterisation of the travelling wave speed in the incompressible limit $\gamma \to \infty$, offering further insight into how the speed of segregated travelling waves is determined jointly by the mobility coefficients of the two cell types, $\mu$ and $\nu$, and the growth kinetics encapsulated by the net growth rate of the density of proliferative cells, $G(p)$.

On the other hand, by revisiting the analysis carried out in~\cite{kim2021interface}, we provide a direct and explicit characterisation of the instability of incompressible segregated circular waves by solving, asymptotically, a shape-perturbation problem for the cellular pressure. Our approach is consistent with the previous analysis, but yields closed-form expressions for the pressure, the interface velocity, and the growth rate of angular modes, thereby making evident the destabilisation mechanisms that may be involved in the formation of finger-like protrusions. We also test our theoretical predictions against numerical experiments.

\subsection{Outline of the paper} The paper is structured as follows. In Sec.~\ref{sec: travelling wave construction}, we recount the construction of one-dimensional segregated travelling waves. In Sec.~\ref{sec: wave speed}, we derive variational estimates for the speed of such travelling waves and obtain the explicit expression of the wave speed in the incompressible limit. In Sec.~\ref{sec: stability}, we analyse stability, formalising a one-dimensional argument as to why segregated travelling waves can be expected to be stable only when $\nu>\mu$, and investigate the stability of incompressible segregated circular waves. In Sec.~\ref{sec:discresperp}, we conclude with a discussion of the results obtained and opportunities for future directions.

\section{The shape of segregated travelling waves}\label{sec: travelling wave construction} 
In this section, we follow the approach in~\cite{lorenzi2016interfaces} and summarise the construction of segregated travelling wave solutions for Eqs.~\eqref{eq:full_model} in one spatial dimension, \textit{i.e.}, $\mathbf{x} \equiv x \in \mathbb{R}$. We first consider a general pressure law that satisfies the assumptions given by Eq.~\eqref{ass:pressure} (see Sec.~\ref{sec: travelling wave constructiongp}). Next, we turn to the case of the power-law defined via Eq.~\eqref{eq: pl pressure} (see Sec.~\ref{sec: travelling wave constructionplp}). We then formally investigate the incompressible regime (see Sec.~\ref{sec: travelling wave constructionplpincomp}).

\subsection{The shape of segregated travelling wave solutions for a general pressure law}\label{sec: travelling wave constructiongp} 
Under the assumptions on the net growth rate $G$, given by Eq.~\eqref{ass:G}, and the assumptions on the pressure law $\Pi$, given by Eq.~\eqref{ass:pressure}, we seek one-dimensional travelling wave solutions with constant speed $c>0$. We look for travelling wave solutions with strictly segregated components such that the non-proliferative cells are confined to a finite region at the leading edge, where the pressure approaches zero, while the proliferative cells are restricted to the rear of the wave, where the pressure saturates at the homeostatic value $\bar p$. These solutions are of the form: 
 \begin{equation}
 \label{eq:TWansatz}
(\rho(x,t), \eta(x,t)) = (\uu(z),\vv(z)) \, , \quad z=x-ct\,,
\end{equation} 
such that, up to a translation and for some $\ell > 0$ to be determined,
 \begin{equation}
 \label{eq:SegrTW}
\begin{split}
&\uu(z) > 0, \quad \vv(z) = 0 \quad \text{for} \, z \in (-\infty, 0)\,,\\
&\uu(z) = 0, \quad \vv(z) > 0 \quad \text{for} \,z \in [0, \ell) \, ,\\
&\uu(z) = 0, \quad \vv(z) = 0 \quad \text{for} \,z \in [\ell, \infty) \, ,\\
\end{split}
\end{equation} 
and satisfy the boundary conditions
\begin{equation}
\label{eq:BCsTW}
\uu(-\infty) = \Pi^{-1}(\bar p) \, , \quad \vv(\ell) = 0 \, ,
\end{equation} 
where $\Pi^{-1}$ is the inverse\footnote{Note that $\Pi^{-1}$ is well defined as the monotonicity assumption given by Eq.~\eqref{ass:pressure} ensures that the function $\Pi$ is invertible.} of $\Pi$. In addition, the complementary condition
\begin{equation}
\label{eq:compCTW}
\dfrac{\dd\uu}{\dd z}(-\infty) = 0 \, ,
\end{equation} 
is imposed, corresponding to the fact that the profile of $\uu(z)$ flattens as it approaches the equilibrium state $\uu=\Pi^{-1}(\bar p)$. Note that, under the segregation conditions in Eq.~\eqref{eq:SegrTW}, conservation of mass ensures that
\begin{equation}
\label{eq:massconsTW}
    \int_0^\ell \vv(z) \,\mathrm{d}z = m_\vv\,,
\end{equation} 
where the parameter $m_\vv>0$ represents the total mass of the non-proliferative cells.

Substituting the travelling wave ansatz given by Eq.~\eqref{eq:TWansatz} into Eqs.~\eqref{eq:full_model} and imposing the segregation conditions given by Eq.~\eqref{eq:SegrTW} produces the following system:
\begin{subequations}\label{eq:full_model_TW}
\begin{align}
   c\dfrac{\dd\uu}{\dd z}& +\mu \dfrac{\dd}{\dd z}\left(\uu\dfrac{\dd p}{\dd z}\right) + \uu G(p)=0\,, \quad p := \Pi(\uu)\, ,  \quad z \in (-\infty, 0) \, ,\label{eq:eq_u_TW}
    \\[0.2cm]
    c\dfrac{\dd\vv}{\dd z}& +\nu \dfrac{\dd}{\dd z}\left(\vv\dfrac{\dd p}{\dd z}\right) =0\,, \quad p := \Pi(\vv)\, , \quad z \in [0, \ell) \, .\label{eq:eq_v_TW}
\end{align}
\end{subequations}
We complement Eqs.~\eqref{eq:full_model_TW} with the conditions given by Eqs.~\eqref{eq:BCsTW} and~\eqref{eq:compCTW}. For fixed $c>0$, integrating Eq.~\eqref{eq:eq_v_TW} over $(0, \ell)$ and imposing the boundary condition $\vv(\ell)=0$ yields 
\begin{equation} \label{eq:poell}
c + \nu \dfrac{\dd p}{\dd z}=0\,, \quad z \in [0, \ell) \, ,
\end{equation}
implying that the pressure gradient satisfies 
\begin{equation*} 
     \frac{\mathrm{d}p(0^+)}{\mathrm{d}z} =  -\frac{c}{\nu}\,. 
\end{equation*}
Solving Eq.~\eqref{eq:poell} subject to the boundary condition $p(\ell)=0$ (cf. the boundary condition on $\vv$ given by Eq.~\eqref{eq:BCsTW}), and recalling the segregation conditions in Eq.~\eqref{eq:SegrTW}, which require that $\vv \equiv 0$ on $(-\infty,0) \cup [\ell, \infty)$, yields:
\begin{equation}\label{eq:eta_explicit}
    \vv(z) = \begin{cases}
         \Pi^{-1}\left(\dfrac{c}{\nu}(\ell-
    z)\right)\,,\quad&\mbox{for } z\in[0,\ell]\,,
        \\[0.2cm]
        0\,,\quad & \mbox{for } z\in (-\infty,0) \cup (\ell, \infty) \, .
\end{cases}
\end{equation}
The value of $\ell$ is determined through Eq.~\eqref{eq:massconsTW}, that is,
\begin{equation}
\label{eq:ellgenericPi}
    \int_0^\ell \Pi^{-1}\left(\frac{c}{\nu}(\ell-z)\right)\,\mathrm{d}z = m_\vv\,.
\end{equation}

We then turn to the differential equation for the density of proliferative cells $\uu$, that is Eq.~\eqref{eq:eq_u_TW}. For $c>0$, under the assumptions given by Eqs.~\eqref{ass:G} and~\eqref{ass:pressure}, the solutions of Eq.~\eqref{eq:eq_u_TW} subject to the conditions given by Eqs.~\eqref{eq:BCsTW} and~\eqref{eq:compCTW} are positive, continuous and, by the maximum principle, monotonically decreasing on $(-\infty,0)$. Moreover, from~\cite{lorenzi2016interfaces} we know that $p$ is such that
\begin{equation}
\label{eq:propinterfp}
p(0^-) = p(0^+) = \dfrac{c}{\nu} \ell \, , \quad \frac{\mathrm{d}p(0^-)}{\mathrm{d}z} = \frac{\nu}{\mu} \frac{\mathrm{d}p(0^+)}{\mathrm{d}z} = -\frac{c}{\mu} \, ,
\end{equation}
and, therefore, 
$$
\uu(0^-) = \Pi^{-1}\left(\dfrac{c}{\nu} \ell \right) \, , \quad \frac{\mathrm{d}\uu(0^-)}{\mathrm{d}z} = -\frac{c}{\mu} \, \dfrac{1}{\Pi'(\uu(0^-))} \, .
$$
Hence, for a fixed $c>0$, $\uu(z)$ is a monotonically decreasing solution to the following problem
\begin{subequations}\label{eq:spec_rho_eq_z}
\begin{align}
     & \mu \dfrac{\dd}{\dd z}\left(\uu\dfrac{\dd p}{\dd z}\right)  + c\dfrac{\dd\uu}{\dd z} + \uu G(p)=0\,, \quad \text{for}\,\,\, z\in(-\infty,0)\,,\\[0.2cm] 
     & \uu(0) = \Pi^{-1}\left(\dfrac{c}{\nu} \ell \right) \, , \quad \frac{\mathrm{d}\uu(0)}{\mathrm{d}z} = -\frac{c}{\mu} \, \dfrac{1}{\Pi'(\uu(0))} \, ,
     \end{align}
\end{subequations}
where $p := \Pi(\uu)$.
Notably, while the aforementioned construction of segregated travelling waves is formally independent of the relative value of the mobility coefficients $\mu$ and $\nu$, the relative value of these parameters impacts the stability of such waves, as later demonstrated in Sec.~\ref{sec: stability}.

\subsection{The shape of segregated travelling waves for the pressure law defined via Eq.~\eqref{eq: pl pressure}}
\label{sec: travelling wave constructionplp}
Substituting Eq.~\eqref{eq: pl pressure} into Eqs.~\eqref{eq:eta_explicit}-\eqref{eq:spec_rho_eq_z}, renaming $\uu$, $\vv$, $p$, $c$, and $\ell$ to $\uu_\gamma$, $\vv_\gamma$, $p_\gamma$, $c_\gamma$, and $\ell_\gamma$, respectively, in order to highlight dependence on the pressure law parameter $\gamma$, and introducing the notation
\begin{equation}\label{def:rhocgamma}
\varrho_{c_\gamma,\gamma} := \left(\dfrac{m_{\vv_\gamma} \, c_{\gamma}}{\nu}\right)^\frac{1}{\gamma+1} \, ,
\end{equation}
we find that
\begin{equation}\label{eq:eta_explicit_powerlaw}
    \vv_{\gamma}(z) = \begin{cases}
\varrho_{c_{\gamma},\gamma} \left(1-\dfrac{z}{\ell_\gamma}\right)^\frac{1}{\gamma}\,,\quad&\mbox{for } z\in[0,\ell_\gamma)\,, 
        \\[0.2cm]
        0\,,\quad & \mbox{for } z\in (-\infty,0) \cup [\ell_\gamma, \infty) \, ,
\end{cases}
\end{equation}
\begin{equation}\label{eq:ell_explicit_powerlaw}
m_{\vv_\gamma}=\left(\dfrac{\gamma}{\gamma+1}\right)^{\frac{\gamma + 1}{\gamma}} \ell_\gamma^{\frac{\gamma+1}{\gamma}}\left(\dfrac{c_{\gamma}}{\nu}\right)^\frac{1}{\gamma} \quad \Longrightarrow \quad \ell_\gamma = \dfrac{\gamma+1}{\gamma} \, 
\left(\dfrac{\nu}{c_{\gamma}}\right)^{\frac{1}{\gamma+1}} m_{\vv_\gamma}^{\frac{\gamma}{\gamma+1}} \, ,
\end{equation}
and
\begin{subequations}\label{eq:spec_rho_eq_z_powerlaw}
\begin{align}
     & \mu \dfrac{\dd}{\dd z}\left(\uu_{\gamma}\dfrac{\dd p_{\gamma}}{\dd z}\right)  + c_{\gamma}\dfrac{\dd\uu_{\gamma}}{\dd z} + \uu_{\gamma} G(p_{\gamma})=0\,, \quad z\in(-\infty,0)\,,\\[0.2cm] 
     & \uu_{\gamma}(0) = \varrho_{c_{\gamma},\gamma} \, , \quad \frac{\mathrm{d}\uu_{\gamma}(0)}{\mathrm{d}z} = -\frac{c_{\gamma}}{\mu(\gamma+1)\varrho_{c_{\gamma},\gamma}^{\gamma-1}} \, ,
     \end{align}
\end{subequations}
where $p_{\gamma} := \Pi(\uu_{\gamma})$, with $\Pi(\uu_{\gamma})$ defined via Eq.~\eqref{eq: pl pressure}. Eq.~\eqref{eq:spec_rho_eq_z_powerlaw} sets the stage for the analysis of the wave speed $c_{\gamma}$, as shown in Sec.~\ref{sec: wave speed}.

\subsection{The shape of segregated travelling waves in the incompressible limit}
\label{sec: travelling wave constructionplpincomp}
In the vein of~\cite{perthame2014hele}, we now investigate the incompressible limit by formally letting $\gamma \to \infty$. In this asymptotic regime, the pressure saturates, and the problem reduces to a sharp-interface free boundary problem, where the travelling wave structure is fully determined by the interface conditions. In fact, letting $\gamma\to\infty$ in $p_{\gamma} := \Pi(\sigma_\gamma)$ with $\Pi(\sigma_\gamma)$ defined via Eq.~\eqref{eq: pl pressure}, one formally finds that the limiting pressure, $p_\infty$, and the limiting total density, $\sigma_\infty := \uu_\infty + \vv_\infty$, satisfy the Hele--Shaw type constraint
\begin{equation*}
    p_\infty(1-\sigma_\infty)=0 \, .
\end{equation*}
A constraint of this type expresses the fact that either the medium is saturated, namely $\sigma_\infty=1$, or the pressure vanishes, corresponding to the emergence of a sharp-interface configuration separating fully saturated and empty regions. 

As a consequence, under the power-law given by Eq.~\eqref{eq: pl pressure} in the asymptotic regime $\gamma \to \infty$, the components of travelling wave solutions $(\uu_\gamma,\vv_\gamma)$ that satisfy the segregation properties given by Eq.~\eqref{eq:SegrTW} take the form
\begin{equation*}
    \uu_\infty(z)=
    \begin{cases}
1 \,,\quad&\mbox{for } z\in(-\infty,0)\,,
        \\[0.2cm]
        0\,,\quad & \mbox{for } z\in [0, \infty) \, ,
    \end{cases}
    \qquad
    \vv_\infty(z)= \begin{cases}
1 \,,\quad&\mbox{for } z\in[0,\ell_\infty)\,,
        \\[0.2cm]
        0\,,\quad & \mbox{for } z\in (-\infty,0) \cup [\ell_\infty, \infty) \, .
\end{cases}
\end{equation*}
Note that the above expression of $\vv_\infty(z)$ is consistent with the one obtained by formally letting $\gamma \to \infty$ in Eq.~\eqref{eq:eta_explicit_powerlaw}, having noted that letting $\gamma \to \infty$ in Eq.~\eqref{def:rhocgamma} formally gives $\varrho_{c,\infty} = 1$. Moreover, letting $\gamma \to \infty$ in Eq.~\eqref{eq:ell_explicit_powerlaw} yields $\ell_\infty = m_{\vv_\infty}$.
Furthermore, in the asymptotic regime $\gamma \to \infty$, Eq.~\eqref{eq:poell} formally gives  
$$
c_{\infty} + \nu \dfrac{\dd p_{\infty}}{\dd z}=0\,, \quad z \in [0, \ell_\infty) \, .
$$
Solving the above differential equation subject to the boundary condition $p_{\infty}(\ell_\infty)=0$, since $\vv_\infty(\ell)=0$, yields
\begin{equation*}
    p_\infty(z)=\frac{c_{\infty}}{\nu}(\ell_{\infty}-z),\qquad z\in[0,\ell_{\infty}) \, .
\end{equation*}
Therefore, in the asymptotic regime $\gamma \to \infty$, Eq.~\eqref{eq:spec_rho_eq_z_powerlaw} written in terms of the cellular pressure can be formally reduced to 
\begin{subequations}\label{eq:incompressible_pressure}
\begin{align}
     &  \mu \frac{\dd^2 p_\infty}{\dd z^2} + G(p_\infty)=0\,, \quad z\in(-\infty,0)\,,\label{eq:incompressible_pressurea} \\[0.2cm] 
     & p_\infty(0) = \dfrac{c}{\nu} \ell_\infty  \, , \quad \frac{\mathrm{d}p_\infty(0)}{\mathrm{d}z} = -\frac{c_\infty}{\mu} \, .
     \end{align}
\end{subequations}

\section{The speed of segregated travelling waves}\label{sec: wave speed}
In this section, we characterise the invasion speed of segregated travelling waves using variational estimates. We then consider the incompressible limit $\gamma \to \infty$, where the estimates become sharp and provide an explicit formula for the speed of the travelling wave, which is formally shown in Sec.~\ref{sec:wsinclim}. Importantly, the derivation of such variational estimates exploits the fact that Eq.~\eqref{eq:spec_rho_eq_z_powerlaw} can be recast as a travelling wave problem for a generalised porous--Fisher equation with a moving boundary, for which a variational principle providing almost explicit bounds on the wave speed can be derived, as we begin by showing in Sec.~\ref{sec:var_princ_1D}. 

\subsection{Variational wave speed estimates for a generalised porous--Fisher equation with a moving boundary}\label{sec:var_princ_1D}
The porous--Fisher model with a moving boundary was initially introduced and studied by Fadai and Simpson in~\cite{fadai2020new}, whilst a modified version that includes a general diffusivity $D(\rho)$, growth term $f(\rho)$, and boundary density $\varrho$ is instead considered here. In one spatial dimension, the model is given by the following free boundary problem for the (normalised) cell density $\rho(x,t)$:
\begin{subequations}
    \begin{align}
    & \partial_t \rho = \partial_x\left(D(\rho)\,\partial_x\rho\right)+f(\rho)\,, \quad x < s(t)\,,\\
    & \rho = \varrho \, ,\quad \frac{\mathrm{d}s}{\mathrm{d}t} = -\kappa D(\rho)\partial_x \rho\,,\quad x = s(t)\,,
\end{align}\label{eq:fisher_stefan}
\end{subequations}
where $\kappa\geq 0$ is a parameter that controls the speed of the moving boundary $s(t)$. We let the growth term satisfy
\begin{equation}\label{ass:f}
f(0) = f(1) = 0 \, , \quad f(\rho) > 0 \;\; \text{for } \, \rho\in (0,1) \, ,
\end{equation}
and we make the following assumptions on the diffusivity 
$$
D(0) = 0 \, , \quad D(\rho) > 0 \;\; \text{for } \, \rho\in (0,1] \, .
$$
Furthermore, we assume that the boundary density satisfies $0 \leq \varrho < 1$.

Introducing the travelling wave coordinate $z = x-s(t) = x-ct$, where $c\geq0$ is the constant speed of the invading wave, having fixed the moving boundary at $z=0$, and substituting the travelling wave ansatz  $\rho(x,t)=\uu(z)$ into the above free boundary problem defined by Eq.~\eqref{eq:fisher_stefan} yields 
\begin{subequations}
\begin{align}
&\frac{\mathrm{d}}{\mathrm{d}z}\!\left(D(\uu)\frac{\mathrm{d}\uu}{\mathrm{d}z}\right) + c\frac{\mathrm{d}\uu}{\mathrm{d}z} + f(\uu) = 0\, , \quad z\in(-\infty,0)\,, \\[0.2cm] 
&\;\uu(0) = \varrho \, , \quad \frac{\mathrm{d}\uu}{\mathrm{d}z}(0) = -\frac{c}{\kappa D(\varrho)} \, .
\end{align}\label{f-s-tw}
\end{subequations}
Similarly to the methods in~\cite{fadai2020new}, we consider monotone solutions of Eq.~\eqref{f-s-tw} that decrease from $1$ to $\varrho$. We then define $\ww(\uu):=-\mathrm{d}\uu/\mathrm{d}z$ and rewrite Eq.~\eqref{f-s-tw} as 
\begin{subequations}\label{eq:f-s-tw}
    \begin{align}
   & \ww(\uu)\frac{\mathrm{d}}{\mathrm{d}\uu}(D(\uu) \ww(\uu))  - c\ww(\uu) + f(\uu)  = 0\,, \quad \uu \in (\varrho,1],  \label{f-s-tw-va}
    \\[0.2cm]
   & \ww(\varrho) = \frac{c}{\kappa D(\varrho)}\,.
\end{align}
\end{subequations}

Following the ideas in \cite{benguria1994validity,benguria1996speed,stokes2024speed,benguria1996variational, benguria2004minimal, crossley2026optimalcontrolapproachnonlinear}, to obtain a variational characterisation of the travelling wave speed $c$,  we introduce a non-negative, non-increasing test function $\varphi(\uu)$ with $\phi(\uu)=-\varphi'(\uu)\geq0$, chosen such that the endpoint contribution at $\uu=1$ vanishes. 
Multiplying Eq.~\eqref{f-s-tw-va} by $\varphi(\uu)D(\uu)$, integrating between $\varrho$ and $1$, and using the fact that $\ww(1)=0$ (since the profile of $\uu(z)$ flattens as it approaches the equilibrium state $\uu=1$), we obtain 
\begin{equation}
    \int_{\varrho}^1 \left(-\frac{1}{2}\phi(\uu) \ww^2(\uu)D(\uu)^2+c\ww(\uu)D(\uu)\varphi(\uu)\right)\mathrm{d}\uu  + \frac{c^2 \varphi(\varrho)}{2\kappa^2}= \int_{\varrho}^1 f(\uu) D(\uu)\varphi(\uu)\,\mathrm{d}\uu\,. \label{eq:int}
\end{equation}
The functional on the left-hand side of Eq.~\eqref{eq:int} is a quadratic in $\ww(\uu)$, which can be bounded above to give 
\begin{equation}
   \frac{c^2}{2} \left[\frac{\varphi(\varrho)}{\kappa^{2}}+\int_{\varrho}^1 \frac{\varphi^2(\uu)}{\phi(\uu)}\,\mathrm{d}\uu\right] \geq\int_{\varrho}^1 f(\uu) D(\uu)\varphi(\uu)\,\mathrm{d}\uu \,,\label{eq:bound}
\end{equation}
with equality when
\begin{equation*}
    \ww(\uu) = \frac{c \, \varphi(\uu)}{D(\uu)\phi(\uu)}\,,\end{equation*}
    or equivalently when
    \begin{equation*}\varphi(\uu) \propto \exp\left\{-\int_{\varrho}^\uu\frac{c}{\ww(q)D(q)} \, \mathrm{d}q\right\}.
\end{equation*}
The moving boundary also provides the compatibility condition
\begin{equation*}
    \varphi(\varrho) = \phi(\varrho)/\kappa\,. 
\end{equation*}
The simplest test function\footnote{This type of test function produces excellent estimates for the wave speed in the classical Fisher--Stefan problem analysed in \cite{crossley2026optimalcontrolapproachnonlinear}.} we could employ that satisfies these bounds, along with the compatibility condition, is 
\begin{equation}
\label{def:varphiu}
\varphi(\uu) := e^{-\kappa (\uu-\varrho)} \, ,
\end{equation}
which gives
\begin{equation}
   \int_{\varrho}^1 \frac{\varphi^2(\uu)}{\phi(\uu)}\,\mathrm{d}\uu = \frac{1}{\kappa^2}\left(1-e^{-\kappa(1-\varrho)}\right)\,. \label{eq:phi_int}
\end{equation}
The integral on the right-hand side of Eq.~\eqref{eq:bound} requires a specific choice of $D(\uu)$ and $f(\uu)$. Taking 
\begin{equation}
\label{def:DPF}
D(\uu):=\uu^m \, , \quad \text{with} \;m\ge0 \, ,
\end{equation}
we can derive different estimates for the travelling wave speed depending on the choice of the growth term $f(\uu)$. 

As a first example, we consider a logistic-type growth term of the form 
\begin{equation}
\label{def:fthetaPF}
f(\uu) :=\uu(1-\uu^\theta) \, , \quad \text{with} \;\theta>0 \, ,
\end{equation}
that is commonly referred to as the $\theta$-logistic or Richards' growth model \cite{Simpson2022ParameterIdentifiability, Liu2024ParameterIdentifiabilityPDE}.
Eq.~\eqref{def:fthetaPF} provides a nonlinear generalisation of the classical logistic growth model, which can be recovered for $\theta=1$. The $\theta$-logistic growth model is widely used across ecology and in the study of tumour growth to describe non-standard saturation mechanisms and nonlinear crowding effects. The parameter $\theta$ controls the strength of the density-dependent inhibition, where values of $\theta>1$ correspond to a sharper saturation near the carrying capacity (here normalised to 1), whereas $0<\theta<1$ describes an earlier and smoother onset of saturation effects. When $D(\uu)$ is defined via Eq.~\eqref{def:DPF} and $f(\uu)$ is defined via Eq.~\eqref{def:fthetaPF}, the right-hand side of Eq.~\eqref{eq:bound} becomes
\begin{equation*}
\begin{split}
 \int_{\varrho}^1 f(\uu) D(\uu)\varphi(\uu)\,\mathrm{d}\uu=\dfrac{e^{\kappa\varrho}}{\kappa^{m+2}}\left[\Gamma(m+2,\kappa\varrho)-\Gamma(m+2,\kappa)-\dfrac{1}{\kappa^\theta}\Gamma(m+2+\theta,\kappa\varrho)+\dfrac{1}{\kappa^\theta}\Gamma(m+2+\theta,\kappa)\right] \, .  
 \end{split}
\end{equation*}
The resulting estimate for the wave speed is
\begin{equation}\label{c_logistic}
    c^2\ge \dfrac{e^{\kappa\varrho}}{\kappa^m(1-\tfrac{1}{2}e^{\kappa(1-\varrho)})}\left[\Gamma(m+2,\kappa\varrho)-\Gamma(m+2,\kappa)-\dfrac{1}{\kappa^\theta}\Gamma(m+2+\theta,\kappa\varrho)+\dfrac{1}{\kappa^\theta}\Gamma(m+2+\theta,\kappa)\right]\,.
\end{equation}
For small values of $\kappa$, we obtain the expansion
\begin{equation*}
    c^2\ge \dfrac{2\kappa^2}{(m+2)(m+2+\theta)}\left[(m+2+\theta)(1-\varrho^{m+2})-(m+2)(1-\varrho^{m+2+\theta})\right] +O(\kappa^3)\,.
\end{equation*}
When $m = \theta = 1$ and $\varrho = 0$, we obtain $c^2\gtrsim \kappa^2/6 + \mathcal{O}(\kappa^3)$, which is sharp with respect to the asymptotic result by Fadai and Simpson~\cite{fadai2020new}. We observe that the $\theta$-logistic growth term introduces density-dependent inhibition, which decreases the bound for the travelling wave speed. Moreover, larger values of $\theta$ weaken the inhibitory effect, leading to larger wave speeds. 

As a second example, we consider the case in which the growth term is given by the Gompertz law~\cite{Simpson2022ParameterIdentifiability, Liu2024ParameterIdentifiabilityPDE}:
\begin{equation}
\label{def:fgomPF}
f(\uu)=-\uu\log \uu \, .
\end{equation}
Such a growth term is widely used to model biological systems characterised by rapid initial proliferation since density-dependent effects diverge as the density tends to zero. 
As such, the Gompertz growth law has become particularly relevant in tumour growth modelling---see for example~\cite{Alarcon2005MultipleScale}. 
When $D(\uu)$ is defined via Eq.~\eqref{def:DPF} and $f(\uu)$ is defined via Eq.~\eqref{def:fgomPF}, the right-hand side of Eq.~\eqref{eq:bound} can be computed explicitly as
\begin{equation*}
\begin{split}
 \int_{\varrho}^1 f(\uu) D(\uu)\varphi(\uu)\,\mathrm{d}\uu =e^{\kappa\varrho}\dfrac{\partial}{\partial m}\Bigg[\frac{1}{\kappa^{m+2}}\big(\Gamma(m+2,\kappa\varrho)-\Gamma(m+2,\kappa)\big)\Bigg] \,.
 \end{split}
\end{equation*}
Substituting this expression back into Eq.~\eqref{eq:bound} leads to the following estimate for the travelling wave speed:
\begin{equation}\label{c_Gompertz}
c^2\geq \dfrac{\kappa^2 e^{\kappa\varrho}}{1-\tfrac{1}{2}e^{-\kappa(1-\varrho)}}\dfrac{\partial}{\partial m}\Bigg[\frac{1}{\kappa^{m+2}}\big(\Gamma(m+2,\kappa)-\Gamma(m+2,\kappa\varrho)\big)\Bigg]. 
\end{equation}
For small values of $\kappa$, we obtain the expansion
\begin{equation*}
    c^2\ge \dfrac{2\kappa^2}{(m+2)^2}\Big[1-\varrho^{m+1}+(m+2)\varrho^{m+2}\ln\varrho\Big] +O(\kappa^3)\,.
\end{equation*}
Compared to the $\theta$-logistic growth with $\theta=1$ (\textit{i.e.}, classical logistic growth), the Gompertz law yields a larger leading-order small-$\kappa$ wave speed bound. This is consistent with the singular low-density behaviour of the Gompertz per-capita growth rate, since $f(\uu)/\uu=-\log \uu\to\infty$ as $\uu\to0^+$. However, this comparison does not hold uniformly for all values of $\theta$.

\begin{figure}[h!]
    \centering
    \includegraphics[width=\linewidth]{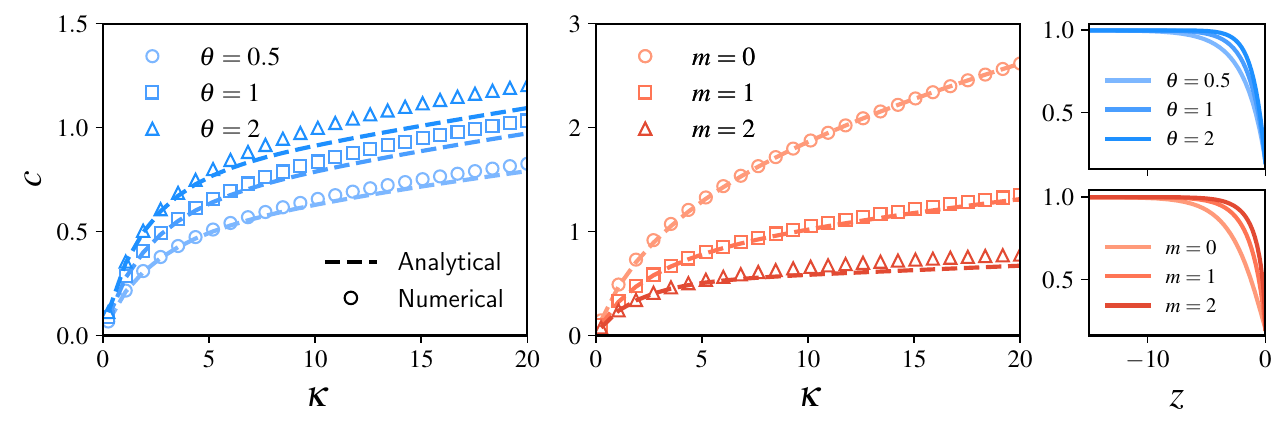}
    \caption{Plots of the travelling wave speeds and solutions for the porous--Fisher model with $\theta$-logistic and Gompertz growth laws. The two larger panels show the numerically estimated travelling wave speed (symbols) and the minimal wave speed predicted by using the variational principle (dashed lines) as a function of $\kappa$. For the $\theta$-logistic growth law in blue (left larger panel), the minimal wave speed is predicted through Eq.~\eqref{c_logistic}; for the Gompertz growth law in red (right larger panel), the prediction of the minimal wave speed is made by using Eq.~\eqref{c_Gompertz}. The two smaller panels display representative travelling wave profiles for each of these growth laws ($\theta$-logistic growth law in blue, top smaller panel, and Gompertz growth law in red, bottom smaller panel) for selected parameter values. Note that here $\varrho=0.2$.}
    \label{fig:pfs-eg}
\end{figure}

To investigate the accuracy of the bounds provided in  Eqs.~\eqref{c_logistic} and~\eqref{c_Gompertz}, numerical simulations are performed and the corresponding wave speeds estimated. To this end, we solve the travelling wave boundary-value problem in Eqs.~\eqref{eq:f-s-tw} using \texttt{solve\_bvp} in \texttt{SciPy}. The wave speed is treated as an additional unknown parameter and computed simultaneously with the wave profile. For each parameter set, continuation in the model parameter is employed by using the previously converged solution as the initial guess.
The resulting plots in Fig.~\ref{fig:pfs-eg} demonstrate that there is a good quantitative agreement between the numerically estimated travelling wave speeds and the minimal wave speeds predicted by Eqs.~\eqref{c_logistic} and~\eqref{c_Gompertz}.

\subsection{Variational estimates for the speed of segregated travelling waves}
\label{sec:var_princ_2species} 
When comparing Eq.~\eqref{eq:spec_rho_eq_z} with Eq.~\eqref{f-s-tw}, renaming $\uu$ and $c$ to $\uu_\gamma$ and $c_\gamma$, we see that the two are equivalent upon setting
\begin{equation}
\label{def:fdrhokappa}
    f(\uu_\gamma) := \uu_\gamma \, G(p_\gamma(\uu_\gamma)) \, ,
\qquad 
D(\uu_\gamma) :=\mu(\gamma+1)\uu^\gamma_\gamma \, ,
\qquad
\varrho:=\varrho_{c_\gamma,\gamma}\, ,
\qquad
\kappa:={\varrho_{c_\gamma,\gamma}^{-1}}\,,
\end{equation}
with $\varrho_{c_\gamma,\gamma}$ defined via Eq.~\eqref{def:rhocgamma}. To ensure consistency between the assumptions given by Eqs.~\eqref{ass:G} and~\eqref{ass:f}, under the pressure law defined via Eq.~\eqref{eq: pl pressure}, we set 
\begin{equation}\label{ass:barp}
\bar p \equiv \bar p_{\gamma} := \frac{\gamma+1}{\gamma}
\end{equation}
so that $G(\bar p) = G(p_\gamma(1)) = 0$.

Choosing the exponential test function $\varphi(\uu_\gamma)$ defined via Eq.~\eqref{def:varphiu}, and substituting the expressions given by Eq.~\eqref{def:fdrhokappa} into Eq.~\eqref{eq:bound}, we obtain
\begin{equation}
 \dfrac{c^2_\gamma}{2} \ge \dfrac{\mu(\gamma+1)}{\varrho_{c_\gamma,\gamma}^2(2-e^{1-1/\varrho_{c,\gamma}})}\int_{\varrho_{c_\gamma,\gamma}}^1\uu_\gamma^{\gamma+1}G(p_\gamma(\uu_\gamma))e^{1-\uu_\gamma/\varrho_{c_\gamma,\gamma}}\,\mathrm{d}\uu_\gamma \, .
 \label{eq: bound two species}
\end{equation}
The estimate in Eq.~\eqref{eq: bound two species} relates the travelling wave speed to the pressure law exponent $\gamma$, the interface density $\varrho_{c_\gamma,\gamma}$, and the growth kinetics encoded in the net growth rate $G$. Different choices of $G$ lead to different lower bounds for the wave speed, thereby allowing for a direct comparison between distinct growth kinetics within the same framework. 

\begin{figure}[htbp]
    \centering
    \includegraphics[width=.8\linewidth]{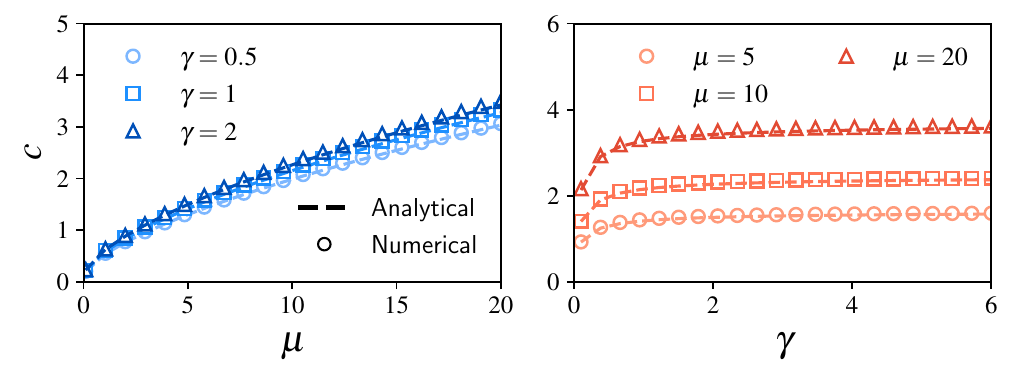}
    \caption{Plots of the speeds of segregated travelling waves from Eq.~\eqref{eq:spec_rho_eq_z} complemented with Eq.~\eqref{def:fthetatwopop}. The numerically estimated travelling wave speed is plotted in symbols and the minimal wave speed predicted using the variational estimate in Eq.~\eqref{eq: logistic two species} is plotted in dashed lines as a function of $\mu$ for different values of $\gamma$ (in blue, left panel) and as a function of $\gamma$ for different values of $\mu$ (in red, right panel). In both panels, $\nu = \mu +1$. As with the generalised porous--Fisher model (see the details in Sec.\ref{sec:var_princ_1D}), the numerical wave speeds were computed by solving the corresponding travelling wave boundary-value problem with continuation in parameter space. 
    }
    \label{fig:placeholder}
\end{figure}

Under the pressure law defined via Eq.~\eqref{eq: pl pressure}, when considering a logistic-type dependence on the cellular pressure that satisfies the assumptions given by Eq.~\eqref{ass:G}, that is, $G(p_\gamma):=1-\frac{p_\gamma}{\bar p}$, defining $\bar p$ via Eq.~\eqref{ass:barp}, we recover the $\gamma$-logistic law 
\begin{equation}
\label{def:fthetatwopop}
    G(p_\gamma(\uu_\gamma)):=1-\uu_\gamma^\gamma \, .
\end{equation}
 Substituting Eq.~\eqref{def:fthetatwopop} into Eq.~\eqref{eq: bound two species}, we obtain a wave speed estimate analogous to that given by Eq.~\eqref{c_logistic}, that is,
\begin{equation}\label{eq: logistic two species}
     c^2_\gamma\geq \dfrac{\mu(\gamma+1)\varrho_{c_\gamma,\gamma}^\gamma}{e^{-1}-\tfrac{1}{2}e^{-1/\varrho_{c_\gamma,\gamma}}} \left\{\,\!\left[\Gamma(\gamma+2,1)-\Gamma\!\left(\gamma+2,{1}/{\varrho_{c_\gamma,\gamma}}\right)\right]
-
\,\varrho_{c_\gamma,\gamma}^{\gamma}\!\left[\Gamma(2\gamma+2,1)-\Gamma\!\left(2\gamma+2,{1}/{\varrho_{c_\gamma,\gamma}}\right)\right]
\right\},
\end{equation}
which can be seen to match very well the numerically estimated travelling wave speed in Fig.~\ref{fig:placeholder}.

Likewise, we can consider the Gompertz law by taking $G(p_\gamma):=-\log (\frac{p_\gamma}{\bar{p}})$ with $\bar p$ defined via Eq.~\eqref{ass:barp}. It follows that 
\begin{equation}
\label{def:fgomtwopop}
G(p_\gamma(\uu_\gamma))=-\gamma \log \uu_\gamma \, 
\end{equation}
and the variational estimate given by Eq.~\eqref{eq: bound two species} becomes
\begin{equation}
   {c_\gamma^2}
\ge
\frac{\mu\gamma(\gamma+1)}
{\varrho_{c_\gamma,\gamma}^2(e^{-1}-\tfrac{1}{2}e^{-1/\varrho_{c_\gamma,\gamma}})}
\left.
\frac{\partial}{\partial s}
\left[
\varrho_{c_\gamma,\gamma}^{s}
\left(
\Gamma\left(s,{1}/{\varrho_{c_\gamma,\gamma}}\right)-\Gamma(s,1)
\right)
\right]
\right|_{s=\gamma+2}. \label{eq: gompertz two species}
\end{equation}

Note also that, in the incompressible limit $\gamma\to \infty$, the right-hand side of Eq.~\eqref{eq: bound two species} simplifies, providing an almost explicit bound on $c_\gamma$. In fact, recalling that $\varrho_{c_\gamma,\gamma}$ is defined via Eq.~\eqref{def:rhocgamma}, we have $\varrho_{c_\gamma,\gamma}-1\sim \frac{1}{\gamma+1}\log({m_{\vv_\gamma} \, c_\gamma}/{\nu})$ as $\gamma \to \infty$. Hence, using the change of variables $\uu_\gamma^\gamma = e^{-s}$ in Eq.~\eqref{eq: bound two species}, we obtain
\begin{align}
 c_\gamma^2 &\ge \dfrac{2\mu(\gamma+1)}{\varrho_{c_\gamma,\gamma}^2(2-e^{1-1/\varrho_{c_\gamma,\gamma}})}\int_{\varrho_{c_\gamma,\gamma}}^1\uu_\gamma^{\gamma+1}G(p_\gamma(\uu_\gamma))e^{1-\uu/\varrho_{c_\gamma,\gamma}}\,\mathrm{d}\uu_\gamma \nonumber
 \\
 &\sim \frac{2\mu(\gamma+1)}{\gamma}\int^{\log({\nu}/{m_{\vv_\gamma} \, c_\gamma})}_0 e^{-s}G(p_\gamma(\uu_\gamma))\,\mathrm{d}s \nonumber
 \\
 & \label{eq: variational incompressible}\sim 2\mu\int_{{m_{\vv_\gamma} \, c_\gamma}/{\nu}}^1 G(p_\gamma)\,\mathrm{d}p_\gamma\,, \quad \text{as } \gamma \to \infty \, .\end{align}
In the next section, we use direct calculations to demonstrate that the bound in Eq.~\eqref{eq: variational incompressible} is sharp and, in fact, predicts the exact wave speed in the incompressible limit. 

\subsection{Exact speed of segregated travelling waves in the incompressible limit} 
\label{sec:wsinclim}
Eq.~\eqref{eq:incompressible_pressure} admits the implicit representation
\begin{equation*}
    \int_{c_\infty\ell_\infty/\nu}^{p_\infty(z)} \frac{\mathrm{d}q}{\sqrt{g(q)}}=-z \, ,\quad\
    g(y):= \frac{2}{\mu} \int_y^{\bar p} G(q)\,\mathrm{d}q \, .
\end{equation*}
Multiplying Eq.~\eqref{eq:incompressible_pressurea} by ${\dd p_\infty}/{\dd z}$ and integrating once yields the first integral
\[
\frac{\mu}{2}\left(\frac{\dd p_\infty}{\dd z}\right)^2
+
\int^{p_\infty}_{\bar p} G(q)\,\mathrm{d}q
=
0,
\]
where we used the asymptotic conditions $p_\infty(-\infty)=\bar p$, ${\dd p_\infty(-\infty)}/{\dd z}=0$ (cf. the conditions given by Eqs.~\eqref{eq:BCsTW} and~\eqref{eq:compCTW}).
Hence,
\[
\int_{c_\infty\ell_\infty/\nu}^{p_\infty(z)}
\frac{\mathrm{d}q}{\sqrt{g(q)}}=-z\, ,
\]
from which, evaluating at the interface with ${\dd p_\infty(0)}/{\dd z}= -c_\infty/\mu$, we obtain the following exact formula for the travelling wave speed in the incompressible limit:
\begin{equation}\label{eq:c_incompressible_general}
    c_\infty^2=2\mu \int_{c_\infty\ell_\infty/\nu}^{\bar p} G(p)\,\mathrm{d}p \, .
\end{equation}
Recalling that $\ell_\infty = m_{\vv_\infty}$ and noting that, when $\bar p \equiv \bar p_{\gamma}$, with $\bar p_{\gamma}$ defined via Eq.~\eqref{ass:barp}, then $\bar p_{\gamma} \to 1$ as $\gamma \to \infty$, Eq.~\eqref{eq:c_incompressible_general} shows that the variational bound given by Eq.~\eqref{eq: variational incompressible} is sharp, thus confirming the validity of the variational principle to predict the speed of segregated travelling waves. The expression given by Eq.~\eqref{eq:c_incompressible_general} also implies uniqueness of the incompressible travelling wave speed, which can be seen via a straightforward monotonicity argument. 

Interestingly, the formula given by Eq.~\eqref{eq:c_incompressible_general} indicates that, in the incompressible regime, the travelling wave speed depends explicitly on the mobility coefficient and total mass of non-proliferative cells through the ratio $\ell_\infty/\nu=m_{\vv_\infty}/\nu$. The expression given by this formula also reveals the existence of an intrinsic upper bound for the wave speed. In fact, in the asymptotic regime of highly mobile non-proliferative cells, namely $\nu\to\infty$, we obtain
\begin{equation*}
       c_\infty^2 \xrightarrow{\nu\to\infty} 2\mu\int_0^{\bar p} G(p)\,\mathrm{d}p=\mathcal{O}(\mu) \, .
\end{equation*}
Hence, even in the limit of arbitrarily large mobility coefficient of non-proliferative cells at the leading edge of the invading front, in the incompressible regime, the propagation speed remains limited by the mobility coefficient of proliferative cells at the rear of the wave and scales as $c_\infty^2\sim\sqrt{\mu}$ when $\nu\to\infty$. 

Finally, we highlight that the above construction is independent of the sign of $\mu-\nu$. Therefore, segregated travelling wave solutions exist regardless of whether the proliferative or non-proliferative cells are more mobile. The issue of stability, however, is separate and requires an independent analysis, as presented in Sec.~\ref{sec: stability}.

As an example, we now compute the wave speed through the formula given by Eq.~\eqref{eq:c_incompressible_general} for a prototypical growth law corresponding to a logistic-type proliferation mechanism, that is, $G(p):= \alpha \, \left(1-\frac{p}{\bar p}\right)$ where $\alpha>0$ is the intrinsic growth rate of the density of proliferative cells. Substituting this definition of $G(p)$ into Eq.~\eqref{eq:c_incompressible_general}, computing the integral and rearranging terms gives the following quadratic equation for $c_{\infty}$:
$$
\left(1 - \dfrac{\alpha \mu}{\bar p} \dfrac{\ell_{\infty}^2}{\nu^2}\right) \, c^2_{\infty} + 2 \alpha \mu \dfrac{\ell_{\infty}}{\nu} \, c_{\infty} - \alpha \mu \bar p = 0 \, .
$$
Solving the above equation and taking the positive root yields
\begin{equation}
\label{eq:c_incompressible_logistic}
    c_{\infty} =\frac{\bar p}{\dfrac{m_{\vv_\infty}}{\nu}+\sqrt{\dfrac{\bar p}{\alpha \mu}}} \, ,
\end{equation}
where we used $\ell_\infty = m_{\vv_\infty}$. 
Eq.~\eqref{eq:c_incompressible_logistic} clearly highlights how the wave speed is determined by the interplay between the different key parameters of the pressure-based model.

\section{Stability of segregated travelling waves}
\label{sec: stability} As mentioned earlier in the paper, the construction of one-dimensional segregated travelling waves is independent of the relative value of the mobility coefficients $\mu$ and $\nu$. However, numerical simulations indicate that, when $\nu < \mu$, segregated travelling waves are unstable in one spatial dimension and finger-like patterns can emerge in two spatial dimensions as a result of the destabilisation of segregated circular waves~\cite{chaplain2020bridging,lorenzi2016interfaces}. To analytically investigate these stability aspects, in Sec.~\ref{sec:stab1d} we formalise an argument as to why one-dimensional segregated travelling waves can be expected to be stable only when $\nu>\mu$. Moreover, in Sec.~\ref{sec:stab2d} we study the stability of incompressible segregated circular waves and shed light on the mechanism of destabilisation that may lead to the emergence of fingering instability.

\subsection{Stability of segregated one-dimensional travelling waves}  \label{sec:stab1d}
A direct way to study the stability of the one-dimensional segregated travelling waves $(\rho(x,t),\eta(x,t)) = (\uu(z),\vv(z))$ considered here is by introducing a small perturbation $\delta\eta(z,t)$ in the density of non-proliferative cells to the left of the interface at $z=0$, while leaving the density of proliferative cells unchanged. In this way, to a first-order approximation, the density of proliferative cells is still $\rho(x,t) = \uu(z)$ while the density of non-proliferative cells now takes the form $\eta(x,t) = \vv(z)$ on $[0,\ell)$ and $\eta(x,t) = \delta\eta(z,t)$ on $(-\infty,0)$. Specifically, we consider perturbations of the form
\begin{equation*}
   \delta \eta(z,t) = \varepsilon e^{\lambda t}\tilde \eta(z)  =\varepsilon e^{\lambda t + z/L}\,, \quad z \in (-\infty,0)
\end{equation*}
for $0<\varepsilon \ll 1$, with $L$ small enough so that the sign of $\lambda$ determines the stability of the wave to perturbations of such form, and $L$ gives the spatial decay of the perturbation moving away from the interface at $z=0$. Substituting into Eq.~\eqref{eq:full_model} for $\eta$, to a first-order approximation in $\varepsilon$, we obtain 
\begin{subequations}
    \begin{align}
\lambda\Tilde{\eta} & = \frac{\dd}{\dd z}\left[\Tilde{\eta}\left(c+\nu\frac{\dd p}{\dd z}\right)\right]\,, \quad z\in(-\infty,0) \, . \nonumber
    \end{align}
\end{subequations}
For $L$ sufficiently small, the perturbation $\tilde{\eta}$ decays on a much shorter length scale than that over which the pressure varies near the interface at $z=0$. Therefore, recalling the properties in Eq.~\eqref{eq:propinterfp}, we have $\dd p/\dd z \sim -c/\mu$ as $z\to 0^-$. The decay rate $\lambda$ can thus be determined by
\begin{equation*}
    \lambda e^{z/L} \sim c\left(1-\frac{\nu}{\mu}\right)\frac{\dd}{\dd z} ( e^{z/L})\,, \quad \text{as } z \to 0^-\,, 
\end{equation*}
which gives
\begin{equation}
    \lambda\sim \frac{c}{L}\left(1-\frac{\nu}{\mu}\right)\, , \quad \text{as } z \to 0^- \, .\label{eq: perturbation}
\end{equation}
This indicates that one-dimensional segregated travelling waves of the form considered here can become unstable when $\nu<\mu$.

\subsection{Stability of incompressible segregated circular waves}  \label{sec:stab2d}
We follow and simplify the stability analysis carried out by Kim and Tong~\cite{kim2021interface} in a nearly radially symmetric setting. Specifically, in two spatial dimensions, whereby $\mathbf{x} \equiv (x,y)$, we consider small perturbations of a circular interface
separating proliferative and non-proliferative cells in the incompressible limit $\gamma \to \infty$, wherein the cell densities are saturated. Without loss of generality, we assume the pressure law $\Pi$ to be such that $\Pi(1)=\bar p$. Hence, denoting the indicator function of the set $\Omega$ by $\chi_\Omega$, letting proliferative cells be surrounded by non-proliferative cells, we take
\[
\rho(x,y,t) =\chi_{\Omega_t} \quad \text{and} \quad \eta(x,y,t)= \chi_{\widetilde{\Omega}_t\setminus\Omega_t},
\]
where $\Omega_t$ is a compact set and $\Omega_t\subset\widetilde{\Omega}_t \subset \mathbb{R}^2$. From now on, for brevity, we omit the subscript denoting time dependence.

We are interested in the dynamics of the interfaces $ \boldsymbol{\gamma}:=\partial\Omega$ (\emph{inner}) and $\widetilde{\boldsymbol{\gamma}} = \partial \widetilde\Omega$ (\emph{outer}), which we parametrise by the polar angle $\theta$. Here $(r,\theta)$ denote the usual polar coordinates. We assume that the interfaces are a small
perturbation of a circle, so that, up to terms of order $\mathcal{O}(\varepsilon)$
for $0<\varepsilon\ll 1$,

\begin{subequations}
\begin{align}
\boldsymbol{\gamma}(\theta,t)
    = R(\theta,t)\,\er \, ,
    \qquad
    R(\theta,t)\sim R_0(t)+\varepsilon R_1(\theta,t) \, ,
\\
\widetilde{\boldsymbol{\gamma}}(\theta,t)
    = \widetilde{R}(\theta,t)\,\er \, ,
    \qquad
    \widetilde{R}(\theta,t)\sim \widetilde{R}_0(t)+\varepsilon \widetilde{R}_1(\theta,t) \,,
\end{align} \label{eq: R asymptotic}
\end{subequations}
where $\er=(\cos\theta,\sin\theta)$. Our goal is to determine the
evolution of the mean radii $R_0(t)$, $\widetilde{R}_0(t)$ and the growth or decay of the perturbations
$R_1(\theta,t)$, $\widetilde{R}_1(\theta,t)$. 

The interfaces dynamics are driven by the gradient of the cellular pressure $p(x,y,t)$, which determines the velocity in the bulk, $\mathbf{u}$, as
\begin{align*}
    \mathbf{u}&=-\mu\nabla p,
    \qquad (x,y) \in\Omega \, ,
\\
\mathbf{u}& =-\nu\nabla p \, ,
    \qquad (x,y) \in\widetilde{\Omega}\setminus\Omega \, .
\end{align*} These expressions are determined by Eqs.~\eqref{eq:full_model} at saturation, which give
\begin{equation}\label{eq: model at saturation}
    -\nabla\cdot\left[\left(\mu \chi_{\Omega} + \nu \chi_{\tilde\Omega\setminus\Omega}\right)\nabla p\right] = \chi_{\Omega} G(p)\,.
\end{equation}
More explicitly, we can write
\begin{subequations}
\begin{align}
    -\mu\Delta p  &= G(p), \quad (x,y)\in{\Omega}\,,
\\
\Delta p &= 0, \quad (x,y)\in\Tilde{\Omega}\setminus\Omega\,.
\end{align}\label{eq: pressure problem}
\end{subequations}
From Eq.~\eqref{eq: model at saturation} we deduce that the pressure and normal velocity are continuous across the inner interface, that is,
\begin{equation}
    p(R^-) = p(R^+)\,,\quad \mu\partial_{\nhat}p(R^-) = \nu\partial_{\nhat} p(R^+)\,,\label{eq: interface conditions}
\end{equation}
where $\nhat$ denotes the outer normal to the inner interface. We also impose that $p(\Tilde{R}) = 0$ for the outer radius.
From Eqs.~\eqref{eq:full_model} we have that the inner interface moves at a normal speed $\vn$ given by
\begin{equation*}
    \vn = \left.\mathbf{u}\right|_{|\mathbf{x}| = R}\cdot \nhat = -\mu\partial_{\nhat} p(R^{-})= -\nu\partial_{\nhat} p(R^{+})\,.
\end{equation*}
This velocity also determines the inner interface dynamics via
\begin{equation*}
    \vn = \partial_t\boldsymbol{\gamma}\cdot\nhat = \partial_t R \,(\er\cdot\nhat)\,,
\end{equation*}
which gives 
\begin{equation}
    \partial_t R\,(\er\cdot\nhat)= -\mu\partial_\nhat p(R^-) = -\nu\partial_\nhat p(R^+)\,.
    \label{eq: Vn R}
\end{equation}

Combining the asymptotic expansion for the inner interface position in Eq.~\eqref{eq: R asymptotic} with the evolution equation in Eq.~\eqref{eq: Vn R}, we seek an asymptotic expansion of the pressure in powers of a small parameter $\varepsilon \ll 1$ that yields an equation for $\partial_t R_1$. We expand 
\begin{equation*}
    p\sim p_0 + \varepsilon p_1\,.
\end{equation*}
Given $\partial_\theta \boldsymbol{\gamma} =\partial_\theta R\er+ R \boldsymbol{\hat e}_\theta$ with $\partial_\theta R\sim\varepsilon \partial_\theta R_1$, the normal vector to the inner interface can be written as
\begin{equation*}
    \nhat = \frac{ R \boldsymbol{\hat e}_r-\partial_\theta R\boldsymbol{\hat e}_\theta}{\sqrt{R^2 + (\partial_\theta R)^2}} = \left[\boldsymbol{\hat e}_r - \varepsilon\frac{\partial_\theta R_1}{R_0}\boldsymbol{\hat e}_\theta + \mathcal{O}(\varepsilon^2)\right]\left[1 + \mathcal{O}(\varepsilon^2)\right]\sim \boldsymbol{\hat e}_r - \varepsilon\frac{\partial_\theta R_1}{R_0}\boldsymbol{\hat e}_\theta \,.
\end{equation*}Then, using the above expansion in Eq.~\eqref{eq: Vn R}, we have $\vn =  \partial_t R + o(\varepsilon)$, or more explicitly
\begin{align}
    \dot R_0 + \varepsilon\partial_t R_1& \sim-\mu\partial_\nhat p_0(R_0^-+\varepsilon R_1^-)-\varepsilon\mu\partial_\nhat p_1(R_0^-+\varepsilon R_1^-)\nonumber
    \\
    &\sim -\mu\partial_\nhat p_0(R_0^-)-\varepsilon\mu\left[R_1\partial_\nhat \partial_r p_0(R_0^-)+\partial_\nhat p_1(R_0^-)\right]\,.\label{eq: asymptotic velocity}
\end{align}
The key observation is that, at order $\varepsilon$, the contribution of $p_1$
to the normal derivative can be evaluated using the radial derivative $\partial_{\nhat}p_1=\partial_r p_1+\mathcal O(\varepsilon)$. 
Note that $p_0$ is radial, as at leading order we have the problem
\begin{align}
-\mu\Delta p_0  = G(p_0), \quad& r \in [0,R_0)\,,\label{eq: p0 inner}
\\
\Delta p_0 = 0, \quad& r \in (R_0, \widetilde{R}_0)\,,\nonumber
\end{align}
with continuous pressure and velocity at the interface, \textit{i.e.}, 
\begin{equation*}
     p_0(R_0^-) = p_0(R_0^+)\,,\quad \mu\partial_rp_0(R_0^-) = \nu\partial_r p_0(R_0^+)\,,
\end{equation*}
and the outer condition $p_0(\Tilde{R}_0) = 0$.
Hence
\begin{align*}
    \partial_\nhat p_0(R_0^-) & = \partial_r p_0(R_0^-)+o(\varepsilon)\,.
\end{align*}
With this, and looking at orders $\mathcal{O}(1)$ and $\mathcal{O}(\varepsilon)$ in Eq.~\eqref{eq: asymptotic velocity}, we obtain 
\begin{subequations}
    \begin{align}
    \dot{R}_0 & = -\mu\partial_r p_0(R_0^-) = -\nu\partial_r p_0 (R_0^+)\,,\label{eq: O1}
    \\
    \partial_t R_1 & = -\mu\partial_r p_1(R_0^-)-\mu R_1\partial_r^2 p_0(R_0^-)\label{eq: Oeps}\,,
\end{align}
\end{subequations}
with the last expression giving the time evolution of $R_1$. The second derivative $\partial_r^2 p_0(R_0^-)$ can be obtained directly from Eq.~\eqref{eq: p0 inner}, writing the Laplacian in polar coordinates and assuming a radial solution, giving
\begin{equation*}
\mu\partial_r^2 p_0 =- G(p_0) - \frac{\mu}{r}\partial_rp_0\,.
\end{equation*}
Evaluating at $r = R_0^-$, and using Eq.~\eqref{eq: O1}, we obtain
\begin{equation}
\mu\partial_r^2 p_0(R_0^-) = -G(p_0(R_0^-)) + \frac{ \dot R_0}{R_0}\,. \label{eq: p0 solution inner}
\end{equation}
We can thus write
\begin{equation}
    \partial_t R_1 = -\mu\partial_r p_1(R_0^-)-R_1\left[-G(p_0(R_0^-)) + \frac{ \dot R_0}{R_0}\right]\,.\label{eq: R1 ev 0}
\end{equation}
The leading order pressure at $r = R_0$ will depend on the outer annulus solution, which at the interface (\textit{i.e.}, for $r\to R_0^+$) gives
\begin{equation}
\nu\partial_r^2 p_0(R_0^+) = \frac{ \dot R_0}{R_0}\,. \label{eq: p0 solution outer}
\end{equation}
For the outer interface we first use conservation of mass for the non-proliferative cells at leading order in $\varepsilon$ to obtain
\begin{equation}
    \widetilde{R}^2_0-{R}^2_0 = \mathcal{A}\,,\quad \dot{\widetilde{R}}_0 = \frac{R_0}{\widetilde{R}_0}\dot{R}_0\,,\label{eq: R0tilde dot}
\end{equation}
for a constant $\mathcal{A}>0$. An analogous argument to the one used for the inner interface yields the following equation for the evolution of $\widetilde{R}_1$:
\begin{equation}
    \partial_t \widetilde{R}_1 = -\nu\partial_r p_1(\widetilde{R}_0^{-})-\frac{\dot{\widetilde{R}}_0}{\widetilde{R}_0}\widetilde{R}_1\,.\label{eq: R1 tilde ev}
\end{equation}

The stability of the interfaces is given by the  order $\mathcal{O}(\varepsilon)$ problem, which determines $\partial_r p_1(R_0^-)$ in Eq.~\eqref{eq: R1 ev 0}, and $\partial_r p_1(\widetilde{R}_0^-)$ in Eq.~\eqref{eq: R1 tilde ev}. Thus, we derive the problem and boundary conditions at order $\mathcal{O}(\varepsilon)$ for $p_1$, which read
\begin{align}
-\mu\Delta p_1  = G'(p_0)p_1, \quad& r \in [0,R_0)\,,\label{eq: p1 inner}
\\
\Delta p_1 = 0, \quad& r \in (R_0, \widetilde{R}_0)\,.\nonumber
\end{align}
The inner interface conditions follow from the asymptotic expansion of Eq.~\eqref{eq: interface conditions}. Pressure continuity along with Eq.~\eqref{eq: O1} give
\begin{equation*}
    p_1(R_0^+) - p_1(R_0^-) = R_1\left[\partial_r p_0(R_0^-)-\partial_rp_0(R_0^+)\right] =R_1\dot R_0\left(\frac{1}{\nu}-\frac{1}{\mu}\right)\,.
\end{equation*}
Similarly, the normal velocity condition gives
\begin{equation*}
    \mu\partial_r p_1(R_0^-)
    -
    \nu\partial_r p_1(R_0^+)
    =
    R_1\left[
        \nu\partial_r^2p_0(R_0^+)
        -
        \mu\partial_r^2p_0(R_0^-)
    \right]\,.
\end{equation*}
Using the radial leading-order solutions on both sides of the interface (cf. Eqs.~\eqref{eq: p0 solution inner} and~\eqref{eq: p0 solution outer}), we obtain
\begin{equation*}
    \mu\partial_r p_1(R_0^-)
    -
    \nu\partial_r p_1(R_0^+)
    =
    R_1G(p_0(R_0^-))\,.
\end{equation*}
Upon expanding the outer interface condition $p(\widetilde{R}) = 0$, we find
\begin{equation*}
    p_1(\widetilde{R}_0) = - \widetilde{R}_1 \partial_r p_0(\widetilde{R}_0^{-}) = \frac{\dot{\widetilde{R}}_0}{\nu}\widetilde{R}_1\,.
\end{equation*}
Summarising, we have found that the perturbation to the interfaces evolves according to
\begin{subequations}\label{eq: R1 time ev}
    \begin{align}
        \partial_t R_1 & = -\mu\partial_r p_1(R_0^-) - R_1\left[-G(p_0(R_0)) + \frac{\dot R_0}{R_0}\right]=-\nu\partial_r p_1(R_0^+) - R_1\frac{\dot R_0}{R_0}\,,\label{eq: R1 ev}
        \\
         \partial_t \widetilde{R}_1 & = -\nu\partial_r p_1(\widetilde{R}_0^{-})-\frac{\dot{\widetilde{R}}_0}{\widetilde{R}_0}\widetilde{R}_1\,,
    \end{align}
\end{subequations}
where
$p_0$ satisfies the leading order (\textit{i.e.}, $\mathcal{O}(1)$) problem
\begin{equation}
    \begin{cases}
-\mu\Delta p_0  = G(p_0), \quad& r \in [0,R_0)\,,
\\
\Delta p_0 = 0, \quad& r \in (R_0, \widetilde{R}_0)\,,
\\
       p_0(R_0^-) = p_0(R_0^+)\,,\\
\mu\partial_rp_0(R_0^-) = \nu\partial_r p_0(R_0^+) = -\dot{R}_0\,,
\\
p_0(\Tilde{R}_0) = 0\,,
\end{cases}\tag{$\mathcal{O}(1)$}
\end{equation}
and $p_1$ solves the first order (\textit{i.e.}, $\mathcal{O}(\varepsilon)$) problem
\begin{equation}
    \begin{cases}
-\mu\Delta p_1  = G'(p_0)p_1, \quad& r \in [0,R_0)\,,
\\
\Delta p_1 = 0, \quad& r \in (R_0, \widetilde{R}_0)\,,
\\p_1(R_0^+) - p_1(R_0^-) =R_1\dot R_0\left[\frac{1}{\nu}-\frac{1}{\mu}\right]
   \,,\\
\nu\partial_r p_1(R_0^+)
    -
    \mu\partial_r p_1(R_0^-)
    =
    -R_1G(p_0(R_0^-))\,,
    \\
    p_1(\Tilde{R}_0) = \dot{\widetilde{R}}_0\widetilde{R}_1/\nu\,,
\end{cases}\tag{$\mathcal{O}(\varepsilon)$}
\end{equation}
with both $p_0$ and $p_1$ being regular at the origin in $r = 0$.

\begin{figure}
    \centering
    \includegraphics[width=0.99\linewidth]{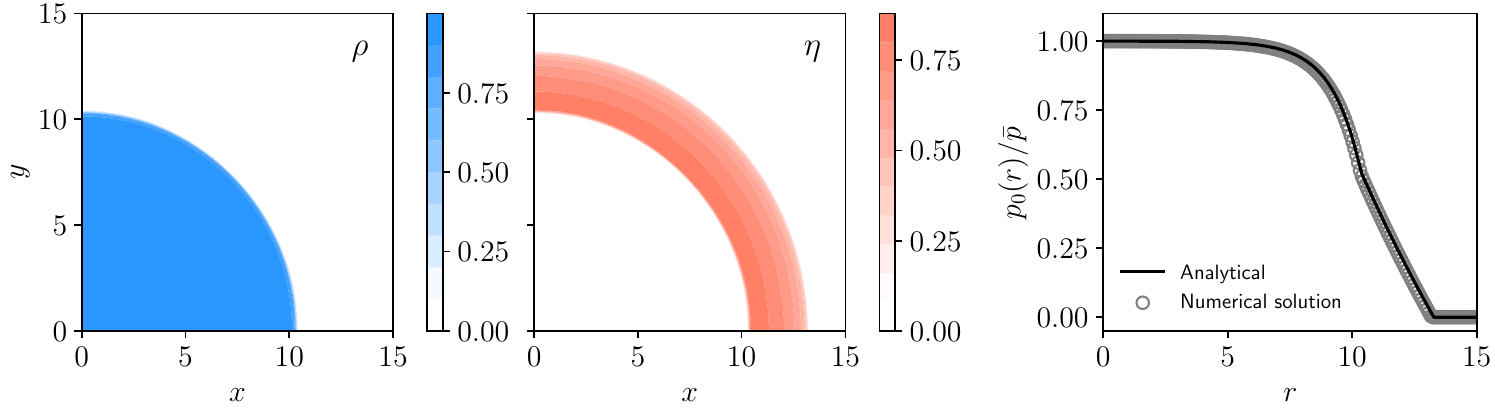}
    \caption{When segregated circular waves propagate, the leading order pressure given by Eq.~\eqref{eq: pressure leading order} approximates the exact profile with the wave speed predicted by Eq.~\eqref{eq: R0dot}. The left panels display the plots of the cell densities $\rho(x,y,t)$ (blue) and $\eta(x,y,t)$ (red) at $t = 10$, and show the propagation of a segregated circular wave. The right panel displays the plot of the leading order pressure $p_0(r)/\overline{p}$. Numerical simulations are carried out setting $\mu = 1$, $\nu = 2$, and $\gamma = 5$, choosing an initial condition corresponding to a segregated circular wave with $R_0(0) = 4$, $\widetilde R_0(0) = 8$, $R_1(0) = 0$, and $\widetilde R_1(0) = 0$, the logistic-type net growth rate $G(p) := 1-p/\bar p$, and using the finite-volume splitting scheme described in~\cite{lorenzi2016interfaces} to solve Eqs.~\eqref{eq:full_model} complemented with the pressure law defined via Eq.~\eqref{eq: pl pressure}.}
    \label{fig: leading order pressure}
\end{figure}

\paragraph{Leading order problem.} The $\mathcal{O}(1)$ problem can be solved in the outer region, where $p_0$ is harmonic, giving
\begin{equation*}
    p_0(r) = -\frac{R_0\dot{R_0}}{\nu}\log\left(\frac{r}{\widetilde{R}_0}\right)\,,
\end{equation*}
where we used that $\nu\partial_rp_0(R_0^+) = -\dot{R}_0$. This gives a closed expression for the pressure at the interface, that is,
\begin{equation*}
    p_0(R_0) = -\frac{R_0\dot R_0}{\nu}\log\left(\frac{R_0}{\widetilde{R}_0}\right)\,.
\end{equation*}
The exact interface speed $\dot{R}_0$ can be evaluated from the inner problem, but in general will depend on the choice of the net growth rate $G$. It is interesting however, as it relates to the one-dimensional problem. We thus solve the problem
\begin{subequations}
    \begin{align}
 -\mu   \partial_r (r\partial_r p_0)   &= rG(p_0),\quad r \in (0,R_0)\,,\label{eq: inner p0}
    \\
    p_0(R_0) & = -\frac{R_0\dot R_0}{\nu}\log\left(\frac{R_0}{\widetilde{R}_0}\right)\,,
    \\
    \mu\partial_r p_0(R_0) & = -\dot R_0\,.
\end{align}
\end{subequations}
We integrate Eq.~\eqref{eq: inner p0} from $0$ to $R_0$ and use regularity at the origin and boundary conditions to find
\begin{equation}
    \dot{R_0} = \frac{1}{R_0}\int_0^{R_0} r G(p_0(r))\,\mathrm{d}r\,.\label{eq: compatibility}
\end{equation}
In general, the equations above cannot be solved explicitly. Analytical progress can be made when the net growth rate $G$ is of the logistic form $G(p) := 1-p/\Bar{p}$, which gives the leading order pressure
\begin{equation}
    \frac{p_0(r)}{\bar p} = 1
-
\frac{\dot R_0}{\sqrt{\mu\bar p}}\,
\frac{
I_0\!\left({r}/\sqrt{\mu\bar p}\right)
}{
I_1\!\left({R_0}/\sqrt{\mu\bar p}\right)
}\,,\label{eq: pressure leading order}
\end{equation}
and the interface radius
\begin{equation}
\label{eq: R0dot}
    \dot R_0
=
\frac{1}{
\dfrac{
I_0\!\left({R_0}/\sqrt{\mu\bar p}\right)
}{
I_1\!\left({R_0}/\sqrt{\mu\bar p}\right)\sqrt{\mu\bar p}
}
+
\dfrac{R_0}{\nu\bar p}
\log\!\left(\dfrac{\widetilde{R}_0}{R_0}\right)
}.
\end{equation}
In particular, this expression shows that the interface speed, which is determined by $\mu$ and $\nu$, is always positive (cf. Eq.~\eqref{eq: compatibility}), so the stability must be determined by the next orders in the pressure. We note that $\widetilde{R}_0$ can be determined as a function of time via Eq.~\eqref{eq: R0tilde dot}. This completely determines the pressure profile (cf. Eq.~\eqref{eq: pressure leading order}) at leading order---see Fig.~\ref{fig: leading order pressure} showing excellent agreement even for moderate values of $\gamma$.

\paragraph{First order problem.} Analytical progress can be made in the $\mathcal{O}(\varepsilon)$ problem by decomposing $R_1$, $\widetilde{R}_1$, and $p_1$ in Fourier modes---note that these are $2\pi$-periodic functions---giving:
\begin{equation*}
     R_1(\theta,t) = \sum_{n\in\mathbb Z}a_n(t)e^{in\theta}\,,\quad \widetilde{R}_1(\theta,t) = \sum_{n\in\mathbb Z}\widetilde{a}_n(t)e^{in\theta}\,,\quad p_1(r,\theta,t) = \sum_{n\in\mathbb Z}b_n(r,t)e^{in\theta}\,.
\end{equation*}
We solve for modes $n\neq 0$, which may drive the fingering instability. Substituting into Eqs.~\eqref{eq: R1 time ev}, we have
\begin{subequations}
    \begin{align}
        \dot a_n & = -\nu\partial_rb_n(R_0^+)- a_n\frac{\dot R_0}{R_0}\,,\label{eq: a_n derivative}
        \\
        \dot{\widetilde{a}}_n & = -\nu\partial_rb_n(\widetilde{R}_0^-)- \widetilde{a}_n\frac{\dot{\widetilde{R}}_0}{\widetilde{R}_0}\,.\label{eq: a_n tilde derivative}
    \end{align}\label{eq: a_n time deriv}
\end{subequations}
Solving for $p_1$ in the outer annulus region (i.e. for $R_0< r< \widetilde{R}_0$), where $\Delta p_1 = 0$, we find
\begin{equation*}
\partial_r^2b_n+\frac{\partial_r b_n}{r}-\frac{n^2b_n}{r^2}=0\,,
\end{equation*}
which, under the boundary condition expanded in Fourier modes $b_n(\widetilde{R}_0) = \dot{\widetilde{R}}_0 \widetilde{a}_n/\nu$, gives
\begin{equation*}
    b_n(r) =B_n\left(
r^{|n|}
-
\frac{\widetilde{R}_0^{2|n|}}{r^{|n|}}
\right)
+
\frac{\dot{\widetilde{R}}_0\widetilde{a}_n}{\nu}
\left(
\frac{\widetilde{R}_0}{r}
\right)^{|n|}\,,
\end{equation*}
for some constant $B_n$. 

In the region occupied by proliferative cells only (\textit{i.e.}, for $0<r<R_0$), the coefficients $b_n$ satisfy
\begin{equation*}
    \partial_r^2b_n+\frac{\partial_r b_n}{r}+\left(\frac{G'(p_0(r))}{\mu} - \frac{n^2}{r^2}\right) b_n = 0\,, 
\end{equation*}
where $p_0$ is determined by the leading order problem and depends on the choice of the net growth rate $G$. In order to close Eq.~\eqref{eq: a_n derivative}, we need $b_n$ near the interface. 

As in the leading-order problem, exact analytical progress can be made when the net growth rate $G$ is of the logistic form $G(p) := 1-p/\Bar{p}$. In this case, the equation for $b_n$ reduces to
\begin{equation*}
    \partial_r^2b_n+\frac{\partial_r b_n}{r}-\left(\frac{1}{\mu\bar p} + \frac{n^2}{r^2}\right) b_n = 0\,,
\end{equation*}
which can be solved in terms of Bessel functions. Imposing regularity at the origin yields the solution
\begin{equation*}
    b_n(r)  = C_n I_{|n|}\left(\frac{r}{\sqrt{\mu\bar p}}\right)\,.
\end{equation*}
The constants $B_n$ and $C_n$ can be determined from the interface conditions, expanded in Fourier modes,
\begin{align*}
    b_n(R_0^+)-b_n(R_0^-) & = a_n \dot R_0\left[\frac{1}{\nu}-\frac{1}{\mu}\right]\,,
    \\
    \nu\partial_r b_n(R_0^+) - \mu\partial_r b_n(R_0^-) & = -a_n G(p_0(R_0))\,.
\end{align*}

\begin{figure}[h!]
    \centering
    \includegraphics[width=0.95\linewidth]{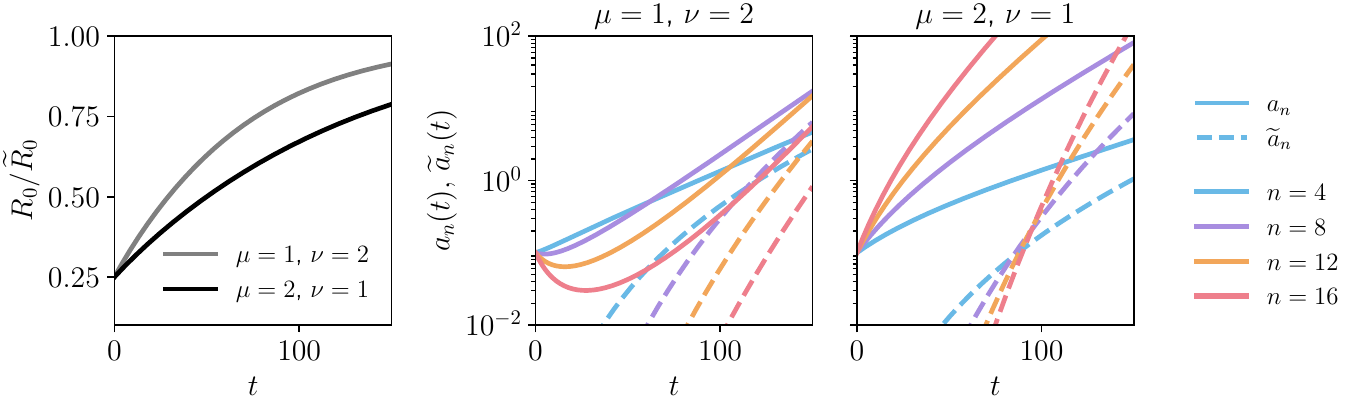}
  \caption{Evolution and stability of segregated circular waves for different relative values of the mobility coefficients. The left panel shows the evolution of the inner-to-outer radius ratio $R_0/\widetilde{R}_0$ for the cases $(\mu,\nu)=(1,2)$ (grey line) and $(\mu,\nu)=(2,1)$ (black line).
In both cases, the ratio approaches 1 as the outer region occupied by non-proliferative cells becomes progressively thinner, although this occurs at different rates in the two cases.
The right panels show the evolution of the perturbation amplitudes
$a_n$ and $\widetilde{a}_n$ for different azimuthal modes $n$: here $a_n(0) = 0.1$ and $\widetilde{a}_n(0) = 0$ for every mode.
For $(\mu,\nu)=(1,2)$, while the corresponding one-dimensional wave is stable, the segregated circular wave eventually becomes unstable: higher-frequency
modes initially decay but subsequently grow as $R_0/\widetilde{R}_0$
approaches 1.
For $(\mu,\nu)=(2,1)$, similarly to the corresponding one-dimensional wave, the segregated circular wave is unstable, with higher-frequency modes growing rapidly and potentially developing into finger-like protrusions at the interface between the two cell types.}
    \label{fig: fig4}
\end{figure}

To solve Eqs.~\eqref{eq: a_n time deriv}, it suffices to determine one of the constants appearing in the expressions for $b_n$ in the two regions. We find the constant $B_n$ associated with the outer solution, which is given by
\begin{equation*}
B_n =
\frac{
-R_0a_nG(p_0(R_0))
+
|n|\dot{\widetilde R}_0\widetilde a_n
\left({\widetilde R_0}/{R_0}\right)^{|n|}
-
\dfrac{\mu}{\nu} R_0\mathcal I_n
\left[
a_n\dot R_0
\left(1-\dfrac{\nu}{\mu}\right)
-
{\dot{\widetilde R}_0\widetilde a_n}
\left({\widetilde R_0}/{R_0}\right)^{|n|}
\right]
}{
\nu |n|
\left(
R_0^{|n|}
+
{\widetilde R_0^{2|n|}}/{R_0^{|n|}}
\right)
-
\mu R_0\mathcal I_n
\left(
R_0^{|n|}
-
{\widetilde R_0^{2|n|}}/{R_0^{|n|}}
\right)
},
\end{equation*}
where
\begin{equation*}
  \mathcal{I}_n = \frac{1}{\sqrt{\mu\bar p}}
\frac{
I_{|n|}'\!\left(
{R_0}/{\sqrt{\mu\bar p}}
\right)
}{
I_{|n|}\!\left(
{R_0}/{\sqrt{\mu\bar p}}
\right)
}>0\,.
\end{equation*}
Next, we use the expression above to rewrite Eqs.~\eqref{eq: a_n time deriv} as the following linear system
\begin{equation}
    \begin{pmatrix}
    \dot{a}_n
    \\
    \dot{\widetilde{a}}_n
    \end{pmatrix}
    = 
    \Lambda(t)\begin{pmatrix}
        a_n
        \\
        \widetilde{a}_n
    \end{pmatrix}\,,\label{eq: 2x2 an system}
\end{equation}
where $\Lambda = (\Lambda_{ij})_{i,j=1}^2$ is a time-dependent matrix obtained by evaluating the left and right partial derivatives of the Fourier coefficients $b_n$ at the inner interface $r = R_0$. The entries of $\Lambda$ read cumbersome, but we provide them in Appendix~\ref{ap: a} for completeness. Although the time-dependent eigenvalues of $\Lambda(t)$ can be calculated, they do not by themselves provide a stability criterion, since the corresponding eigenvectors also change with time. We therefore solve the system in Eq.~\eqref{eq: 2x2 an system} numerically for different parameter values, using Eqs.~\eqref{eq: R0dot} and~\eqref{eq: R0tilde dot} to determine the evolution of $R_0$ and $\widetilde{R}_0$. The dynamics predicted by the system in Eq.~\eqref{eq: 2x2 an system} are illustrated in Fig.~\ref{fig: fig4}. 

Importantly, the stability of the circular wave is not determined solely by the relative value of the mobility coefficients, it also evolves with the geometry of the cell distributions. When $\mu>\nu$, the perturbations grow rapidly, consistently with the
one-dimensional instability criterion. In contrast, when $\mu<\nu$, high-frequency perturbations are
initially damped, but the stabilising effect can be lost as the wave expands and the annular region occupied by the non-proliferative cells becomes thinner,
that is, as $R_0/\widetilde R_0$ approaches 1. Consequently, unlike in one spatial dimension, the condition $\nu>\mu$ alone does not guarantee stability
of a segregated circular wave throughout its evolution.

\begin{figure}[h!]
    \centering
    \includegraphics[width=0.95\linewidth]{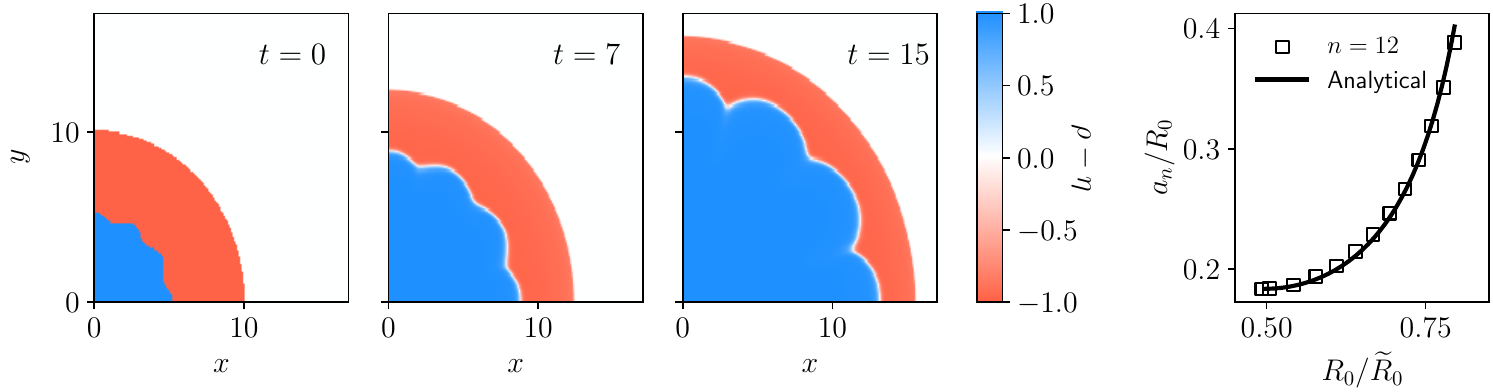}
    \caption{Destabilisation of a segregated circular wave despite the non-proliferative cells being more mobile than the proliferative cells (\textit{i.e.}, for $\nu>\mu$), in contrast with the one-dimensional case. The left panels show the evolution of $\rho(x,y,t)-\eta(x,y,t)$, illustrating the growth of interfacial perturbations and the development of protrusions at the interface between the two cell types. The right panel shows the plot of the amplitude $a_n/R_0$ of the $n=12$ mode as a function of the outer-to-inner radius ratio $R_0/\widetilde{R}_0$. Numerical simulations are carried out setting $\mu = 2$, $\nu = 2.5$, and $\gamma = 25$, choosing an initial condition corresponding to a segregated circular wave with $R_0(0) = 5$, $\widetilde R_0(0) = 10$, $R_1(0) = 0.25$, and $\widetilde R_1(0) = 0$, the logistic-type net growth rate $G(p) := 1-p/\bar p$, and using the finite-volume splitting scheme described in~\cite{lorenzi2016interfaces} to solve Eqs.~\eqref{eq:full_model} complemented with the pressure law defined via Eq.~\eqref{eq: pl pressure}.}
    \label{fig: stable circular wave}
\end{figure}

\begin{figure}
    \centering
    \includegraphics[width=0.95\linewidth]{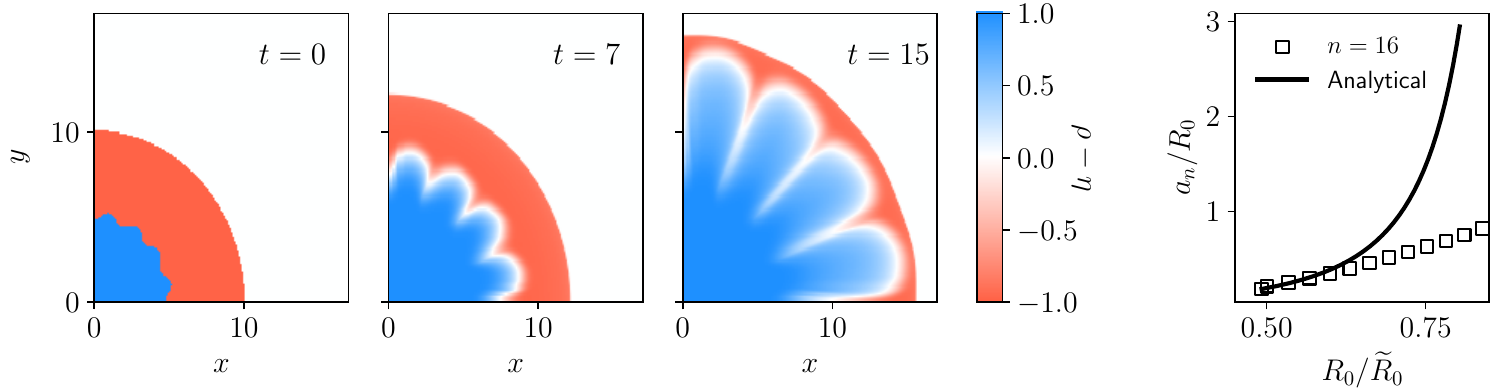}
\caption{Destabilisation of a segregated circular wave when the proliferative cells are more mobile than the non-proliferative cells (\textit{i.e.}, when $\mu>\nu$), consistent also with the instability observed in the one-dimensional case. The left panels show the evolution of $\rho(x,y,t)-\eta(x,y,t)$, illustrating the growth of interfacial perturbations and the development of protrusions at the interface between the two cell types. The right panel shows the plot of the amplitude $a_n/R_0$ of the $n=16$ mode as a function of the outer-to-inner radius ratio $R_0/\widetilde{R}_0$---note that the analytical prediction is an asymptotic result and thus only holds when $a_n\ll R_0$. Numerical simulations are carried out setting $\mu = 2.5$, $\nu = 2$, and $\gamma = 25$, choosing an initial condition corresponding to a segregated circular wave with $R_0(0) = 5$, $\widetilde R_0(0) = 10$, $R_1(0) = 0.25$, and $\widetilde R_1(0) = 0$, the logistic-type net growth rate $G(p) := 1-p/\bar p$, and using the finite-volume splitting scheme described in~\cite{lorenzi2016interfaces} to solve Eqs.~\eqref{eq:full_model} complemented with the pressure law defined via Eq.~\eqref{eq: pl pressure}.}
    \label{fig: unstable circular wave}
\end{figure}

We test these predictions against numerical simulations of the original finite-$\gamma$ system in Eqs.~\eqref{eq:full_model}. The results obtained are summarised by the plots in Figs.~\ref{fig: stable circular wave} and~\ref{fig: unstable circular wave}. The results of numerical simulations in Fig.~\ref{fig: stable circular wave} demonstrate the aforementioned geometry-driven destabilisation in a scenario where $\nu>\mu$: although the
corresponding one-dimensional wave is expected to be stable, perturbations of
the circular interface eventually grow as the outer annulus becomes
sufficiently thin, leading to the formation of small but well-defined protrusions at the interface between the two cell types. The evolution of
the $n=12$ mode is well captured by the asymptotic prediction over the range in which the perturbation remains small. When $\mu>\nu$, as shown in Fig.~\ref{fig: unstable circular wave},
the instability is considerably stronger and protrusions develop rapidly. In this case, the eventual discrepancy between the predicted and numerically measured amplitudes is expected, since the linear shape-perturbation theory ceases to apply once $a_n$ becomes comparable with
$R_0$.

\section{Discussion and perspectives}
\label{sec:discresperp}
In this work, we address open questions concerning the speed and stability of segregated waves in a minimal pressure-based model of heterogeneous cell populations. First, we obtain previously unavailable estimates for the speed of one-dimensional segregated travelling waves, revealing an interesting connection with a generalised porous–Fisher model with a moving boundary \cite{fadai2020new} and recovering an explicit characterisation of the wave speed in the incompressible limit. We then provide a formal one-dimensional stability argument, absent from previous analyses, explaining why stability of segregated waves depends on the relative value of the cell mobility coefficients. Finally, building on numerical observations presented in~\cite{lorenzi2016interfaces} and the stability analysis carried out in~\cite{kim2021interface}, we derive explicit analytical expressions for the growth of interfacial perturbations within incompressible segregated circular waves. Interestingly, these reveal that, in contrast with the one-dimensional case, the stability of such circular waves is not determined solely by the relative value of the mobility coefficients, and thus instabilities may arise irrespective of which cell type has the larger mobility, a behaviour that we confirmed numerically and that, to our knowledge, had not been observed previously.

Several directions remain open for future investigation. From a mathematical perspective, a fully rigorous one-dimensional stability theory for segregated travelling waves remains to be developed beyond the formal argument presented here, while in higher dimensions a nonlinear stability theory is still lacking. In particular, it remains unclear whether the interfacial instabilities identified by our linear analysis persist beyond the nearly radially symmetric configuration and how such instabilities develop nonlinearly into more complex morphologies. It would also be interesting to extend the results presented here in the context of alternative models for pressure-driven cell population dynamics featuring cross-diffusion terms~\cite{bejar2026analysis,campos2025biomechanical}. More fundamentally, pressure-based models of the type considered here are commonly obtained as limiting descriptions in which viscous stresses are neglected. It would therefore be important to understand how retaining non-negligible viscosity modifies the propagation and stability properties we have identified, and to what extent the behaviours obtained in the zero-viscosity regime persist away from this limiting scenario. Such extensions may also reveal whether additional mechanical effects provide a natural regularisation mechanism for interfacial perturbations and influence the selection of emergent modes, which drive the instability of segregated waves and the formation of finger-like protrusions observed in tissue morphogenesis and tumour invasion.

\section*{Acknowledgments}
For the purpose of Open Access, the authors have applied a CC BY public copyright licence to any Author Accepted Manuscript (AAM) version arising from this submission. CF acknowledges support from a Hooke Research Fellowship. RMC acknowledges support from the Engineering and Physical Sciences Research Council (EP/Z534870/1). MC acknowledges support from the European Union - Next Generation EU, Mission 4, Component 1 (CUP: E13C24002380006) and by the National Group of Mathematical Physics (GNFM-INdAM) (CUP: E5324001950001). MC and TL are members of INdAM-GNFM.

\section*{Appendix}

\appendix
\section{$\Lambda(t)$-matrix elements for the stability of incompressible segregated circular waves}\label{ap: a}
We detail here the matrix elements appearing in Eq.~\eqref{eq: 2x2 an system}. First, we denote
\begin{equation*}
    \mathcal{D}_n := \nu |n|
\left(
R_0^{|n|}
+
\frac{\widetilde R_0^{2|n|}}{R_0^{|n|}}
\right)
-
\mu R_0\mathcal I_n
\left(
R_0^{|n|}
-
\frac{\widetilde R_0^{2|n|}}{R_0^{|n|}}
\right)\,.
\end{equation*}Now, from Eq.~\eqref{eq: a_n derivative} we obtain
\begin{align*}
\Lambda_{11}
&=
\frac{\nu |n|}{\mathcal D_n}
\left(
R_0^{|n|}
+
\frac{\widetilde R_0^{2|n|}}{R_0^{|n|}}
\right)
\left[
G(p_0(R_0))
+
\frac{\mu}{\nu}\mathcal I_n\dot R_0
\left(1-\frac{\nu}{\mu}\right)
\right]
-
\frac{\dot R_0}{R_0},
\\[0.5em]
\Lambda_{12}
&=
-\frac{\nu |n|\dot{\widetilde R}_0}{R_0\mathcal D_n}
\left(
R_0^{|n|}
+
\frac{\widetilde R_0^{2|n|}}{R_0^{|n|}}
\right)
\left(\frac{\widetilde R_0}{R_0}\right)^{|n|}
\left(
|n|+\frac{\mu}{\nu}R_0\mathcal I_n
\right) + 
\frac{|n|\dot{\widetilde R}_0}{R_0}
\left(\frac{\widetilde R_0}{R_0}\right)^{|n|}\,,
\end{align*}
while from Eq.~\eqref{eq: a_n tilde derivative} we obtain
\begin{align*}
\Lambda_{21}
&=
\frac{2\nu |n|R_0\widetilde R_0^{|n|-1}}{\mathcal D_n}
\left[
G(p_0(R_0))
+
\frac{\mu}{\nu}\mathcal I_n\dot R_0
\left(1-\frac{\nu}{\mu}\right)
\right],
\\[0.5em]
\Lambda_{22}
&=
-\frac{2\nu |n|\dot{\widetilde R}_0\widetilde R_0^{|n|-1}}
{\mathcal D_n}
\left(\frac{\widetilde R_0}{R_0}\right)^{|n|}
\left(
|n|+\frac{\mu}{\nu}R_0\mathcal I_n
\right)
+
(|n|-1)\frac{\dot{\widetilde R}_0}{\widetilde R_0}\,.
\end{align*}

\bibliography{references}

@article{bejar2026analysis,
	author = {B{\'e}jar-L{\'o}pez, Alexis and Granero-Belinch{\'o}n, Rafael and Pulido, Carlos and Soler, Juan},
	journal = {SIAM Journal on Mathematical Analysis},
	number = {2},
	pages = {2000--2029},
	publisher = {SIAM},
	title = {Analysis of a cross-nonlinear porous-medium system modeling pressure-driven cell population dynamics},
	volume = {58},
	year = {2026}}

@article{shraiman2005mechanical,
	author = {Shraiman, Boris I},
	journal = {Proceedings of the National Academy of Sciences},
	number = {9},
	pages = {3318--3323},
	publisher = {National Academy of Sciences},
	title = {Mechanical feedback as a possible regulator of tissue growth},
	volume = {102},
	year = {2005}}

@article{campos2025biomechanical,
	author = {Campos, Juan and Pulido, Carlos and Soler, Juan},
	journal = {Journal of Nonlinear Science},
	number = {4},
	pages = {87},
	publisher = {Springer},
	title = {Biomechanical effects on traveling waves at the interface of cell populations},
	volume = {35},
	year = {2025}}

@article{meacham2013tumour,
	author = {Meacham, Corbin E and Morrison, Sean J},
	journal = {Nature},
	number = {7467},
	pages = {328--337},
	publisher = {Nature Publishing Group UK London},
	title = {Tumour heterogeneity and cancer cell plasticity},
	volume = {501},
	year = {2013}}

@article{huang2009non,
	author = {Huang, Sui},
	journal = {Development},
	number = {23},
	pages = {3853--3862},
	publisher = {Company of Biologists},
	title = {Non-genetic heterogeneity of cells in development: more than just noise},
	volume = {136},
	year = {2009}}

@article{wang2022cellular,
	author = {Wang, Jingru and He, Jia and Zhu, Meishu and Han, Yan and Yang, Ronghua and Liu, Hongwei and Xu, Xuejuan and Chen, Xiaodong},
	journal = {Stem Cell Reviews and Reports},
	number = {6},
	pages = {1912--1925},
	publisher = {Springer},
	title = {Cellular heterogeneity and plasticity of skin epithelial cells in wound healing and tumorigenesis},
	volume = {18},
	year = {2022}}

@article{rognoni2018skin,
	author = {Rognoni, Emanuel and Watt, Fiona M},
	journal = {Trends in Cell Biology},
	number = {9},
	pages = {709--722},
	publisher = {Elsevier},
	title = {Skin cell heterogeneity in development, wound healing, and cancer},
	volume = {28},
	year = {2018}}

@misc{darcy1856fontaines,
	author = {Darcy, H.},
	howpublished = {V. Dalmont, Paris, 305-311},
	title = {Les Fontaines Publiques de la Ville de Dijon},
	year = {1856}}

@article{wang2012adipose,
	author = {Wang, Yuan-Yuan and Lehu{\'e}d{\'e}, Camille and Laurent, Victor and Dirat, B{\'e}atrice and Dauvillier, St{\'e}phanie and Bochet, Ludivine and Le Gonidec, Sophie and Escourrou, Ghislaine and Valet, Philippe and Muller, Catherine},
	journal = {Cancer Letters},
	number = {2},
	pages = {142--151},
	publisher = {Elsevier},
	title = {Adipose tissue and breast epithelial cells: a dangerous dynamic duo in breast cancer},
	volume = {324},
	year = {2012}}

@article{friedman2015free,
	author = {Friedman, Avner},
	journal = {Philosophical Transactions of the Royal Society A: Mathematical, Physical and Engineering Sciences},
	number = {2050},
	pages = {20140368},
	publisher = {The Royal Society Publishing},
	title = {Free boundary problems in biology},
	volume = {373},
	year = {2015}}

@article{perthame2014hele,
	author = {Perthame, Beno{\^\i}t and Quir{\'o}s, Fernando and V{\'a}zquez, Juan Luis},
	journal = {Archive for Rational Mechanics and Analysis},
	pages = {93--127},
	publisher = {Springer},
	title = {The {H}ele--{S}haw asymptotics for mechanical models of tumor growth},
	volume = {212},
	year = {2014}}

@article{mellet2017hele,
	author = {Mellet, Antoine and Perthame, Beno{\^\i}t and Quir{\'o}s, Fernando},
	journal = {Journal of Functional Analysis},
	number = {10},
	pages = {3061--3093},
	publisher = {Elsevier},
	title = {A {H}ele--{S}haw problem for tumor growth},
	volume = {273},
	year = {2017}}

@article{kim2016free,
	author = {Kim, Inwon C and Perthame, Beno{\^\i}t and Souganidis, Panagiotis E},
	journal = {Nonlinear Analysis},
	pages = {207--228},
	publisher = {Elsevier},
	title = {Free boundary problems for tumor growth: a viscosity solutions approach},
	volume = {138},
	year = {2016}}

@article{david2024incompressible,
	author = {David, Noemi and Schmidtchen, Markus},
	journal = {Communications on Pure and Applied Mathematics},
	number = {5},
	pages = {2613--2650},
	publisher = {Wiley Online Library},
	title = {On the incompressible limit for a tumour growth model incorporating convective effects},
	volume = {77},
	year = {2024}}

@article{bubba2020hele,
	author = {Bubba, Federica and Perthame, Beno{\^\i}t and Pouchol, Camille and Schmidtchen, Markus},
	journal = {Archive for Rational Mechanics and Analysis},
	number = {2},
	pages = {735--766},
	publisher = {Springer},
	title = {Hele-{S}haw limit for a system of two reaction-(cross-) diffusion equations for living tissues},
	volume = {236},
	year = {2020}}

@article{tang2014composite,
	author = {Tang, Min and Vauchelet, Nicolas and Cheddadi, Ibrahim and Vignon-Clementel, Irene and Drasdo, Dirk and Perthame, Beno{\^\i}t},
	journal = {Partial Differential Equations: Theory, Control and Approximation: In Honor of the Scientific Heritage of Jacques-Louis Lions},
	pages = {401--429},
	publisher = {Springer},
	title = {Composite waves for a cell population system modeling tumor growth and invasion},
	year = {2014}}

@article{carrillo2024multipop,
	author = {Carrillo, J.A. and Macfarlane, F.R. and Lorenzi, T.},
	journal = {Bulletin of Mathematical Biology},
	number = {77},
	title = {Spatial segregation across travelling fronts in individual-based and continuum models for the growth of heterogeneous cell populations},
	volume = {87},
	year = {2025}}

@article{chaplain2020bridging,
	author = {Chaplain, Mark AJ and Lorenzi, Tommaso and Macfarlane, Fiona R},
	journal = {Journal of Mathematical Biology},
	pages = {343--371},
	publisher = {Springer},
	title = {Bridging the gap between individual-based and continuum models of growing cell populations},
	volume = {80},
	year = {2020}}

@article{bertsch2015travelling,
	author = {Bertsch, M and Hilhorst, D and Izuhara, H and Mimura, M and Wakasa, T},
	journal = {European Journal of Applied Mathematics},
	number = {3},
	pages = {297--323},
	publisher = {Cambridge University Press},
	title = {Travelling wave solutions of a parabolic-hyperbolic system for contact inhibition of cell-growth},
	volume = {26},
	year = {2015}}

@article{basan2009homeostatic,
	author = {Basan, M and Risler, T and Joanny, J F and Sastre-Garau, X and Prost, J},
	journal = {HFSP Journal},
	number = {4},
	pages = {265--272},
	publisher = {Taylor \& Francis},
	title = {Homeostatic competition drives tumor growth and metastasis nucleation},
	volume = {3},
	year = {2009}}

@article{sherratt2001new,
	author = {Sherratt, Jonathan A and Chaplain, Mark AJ},
	journal = {Journal of Mathematical Biology},
	number = {4},
	pages = {291--312},
	publisher = {Springer},
	title = {A new mathematical model for avascular tumour growth},
	volume = {43},
	year = {2001}}

@article{preziosi2009multiphase,
	author = {Preziosi, Luigi and Tosin, Andrea},
	journal = {Journal of Mathematical Biology},
	pages = {625--656},
	publisher = {Springer},
	title = {Multiphase modelling of tumour growth and extracellular matrix interaction: mathematical tools and applications},
	volume = {58},
	year = {2009}}

@article{lowengrub2009nonlinear,
	author = {Lowengrub, John S and Frieboes, Hermann B and Jin, Fang and Chuang, Yao-Li and Li, Xiangrong and Macklin, Paul and Wise, Steven M and Cristini, Vittorio},
	journal = {Nonlinearity},
	number = {1},
	pages = {R1},
	publisher = {IOP Publishing},
	title = {Nonlinear modelling of cancer: bridging the gap between cells and tumours},
	volume = {23},
	year = {2009}}

@article{ciarletta2011radial,
	author = {Ciarletta, Pasquale and Foret, L and Ben Amar, M},
	journal = {Journal of the Royal Society Interface},
	number = {56},
	pages = {345--368},
	publisher = {The Royal Society},
	title = {The radial growth phase of malignant melanoma: multi-phase modelling, numerical simulations and linear stability analysis},
	volume = {8},
	year = {2011}}

@article{bresch2010computational,
	author = {Bresch, Didier and Colin, Thierry and Grenier, Emmanuel and Ribba, Benjamin and Saut, Olivier},
	journal = {SIAM Journal on Scientific Computing},
	number = {4},
	pages = {2321--2344},
	publisher = {SIAM},
	title = {Computational modeling of solid tumor growth: the avascular stage},
	volume = {32},
	year = {2010}}

@article{ambrosi2002closure,
	author = {Ambrosi, Davide and Preziosi, Luigi},
	journal = {Mathematical Models and Methods in Applied Sciences},
	number = {05},
	pages = {737--754},
	publisher = {World Scientific},
	title = {On the closure of mass balance models for tumor growth},
	volume = {12},
	year = {2002}}

@article{greenspan1976growth,
	author = {Greenspan, H P},
	journal = {Journal of Theoretical Biology},
	number = {1},
	pages = {229--242},
	publisher = {Elsevier},
	title = {On the growth and stability of cell cultures and solid tumors},
	volume = {56},
	year = {1976}}

@article{roose2007mathematical,
	author = {Roose, Tiina and Chapman, S Jonathan and Maini, Philip K},
	journal = {SIAM Review},
	number = {2},
	pages = {179--208},
	publisher = {SIAM},
	title = {Mathematical models of avascular tumor growth},
	volume = {49},
	year = {2007}}

@article{byrne1997free,
	author = {Byrne, HM and Chaplain, Mark AJ},
	journal = {European Journal of Applied Mathematics},
	number = {6},
	pages = {639--658},
	publisher = {Cambridge University Press},
	title = {Free boundary value problems associated with the growth and development of multicellular spheroids},
	volume = {8},
	year = {1997}}

@article{ranft2010fluidization,
	author = {Ranft, Jonas and Basan, Markus and Elgeti, Jens and Joanny, Jean-Fran{\c{c}}ois and Prost, Jacques and J{\"u}licher, Frank},
	journal = {Proceedings of the National Academy of Sciences},
	number = {49},
	pages = {20863--20868},
	publisher = {National Acad Sciences},
	title = {Fluidization of tissues by cell division and apoptosis},
	volume = {107},
	year = {2010}}

@article{drasdo2012modeling,
	author = {Drasdo, Dirk and Hoehme, Stefan},
	journal = {New Journal of Physics},
	number = {5},
	pages = {055025},
	publisher = {IOP Publishing},
	title = {Modeling the impact of granular embedding media, and pulling versus pushing cells on growing cell clones},
	volume = {14},
	year = {2012}}

@article{byrne2009individual,
	author = {Byrne, Helen and Drasdo, Dirk},
	journal = {Journal of Mathematical Biology},
	pages = {657--687},
	publisher = {Springer},
	title = {Individual-based and continuum models of growing cell populations: a comparison},
	volume = {58},
	year = {2009}}

@article{byrne2003modelling,
	author = {Byrne, Helen and Preziosi, Luigi},
	journal = {Mathematical Medicine and Biology: {A} {J}ournal of the IMA},
	number = {4},
	pages = {341--366},
	publisher = {OUP},
	title = {Modelling solid tumour growth using the theory of mixtures},
	volume = {20},
	year = {2003}}

@article{stokes2024speed,
	author = {Stokes, Beth M and Rogers, Tim and James, Richard},
	journal = {Bulletin of Mathematical Biology},
	number = {12},
	pages = {147},
	publisher = {Springer},
	title = {Speed and shape of population fronts with density-dependent diffusion},
	volume = {86},
	year = {2024}}

@misc{crossley2026optimalcontrolapproachnonlinear,
	author = {Rebecca M. Crossley and Carles Falc\'o and Ruth E. Baker},
	title = {An optimal control approach to nonlinear wave speed selection in reaction-diffusion equations}}

@article{lorenzi2016interfaces,
	author = {Lorenzi, Tommaso and Lorz, Alexander and Perthame, Beno{\^\i}t},
	journal = {Kinetic and Related Models},
	number = {1},
	pages = {299--311},
	publisher = {Kinetic and Related Models},
	title = {On interfaces between cell populations with different mobilities},
	volume = {10},
	year = {2016}}

@article{kim2021interface,
	author = {Kim, Inwon C and Tong, Jiajun},
	journal = {Interfaces and Free Boundaries},
	number = {2},
	pages = {191--304},
	title = {Interface dynamics in a two-phase tumor growth model},
	volume = {23},
	year = {2021}}

@article{fadai2020new,
	author = {Fadai, Nabil T and Simpson, Matthew J},
	journal = {Journal of Physics A: Mathematical and Theoretical},
	number = {9},
	pages = {095601},
	publisher = {IOP Publishing},
	title = {{New travelling wave solutions of the porous--Fisher model with a moving boundary}},
	volume = {53},
	year = {2020}}

@article{saffman1958penetration,
	author = {Saffman, P. G. and Taylor, G. I.},
	doi = {10.1098/rspa.1958.0085},
	journal = {Proceedings of the Royal Society A},
	number = {1242},
	pages = {312--329},
	title = {{The penetration of a fluid into a porous medium or Hele-Shaw cell containing a more viscous liquid}},
	volume = {245},
	year = {1958}}

@article{wang2021budding,
	author = {Wang, Shaohe and Matsumoto, Kazue and Lish, Samantha R. and Cartagena-Rivera, Alexander X. and Yamada, Kenneth M.},
	doi = {10.1016/j.cell.2021.05.015},
	journal = {Cell},
	number = {14},
	pages = {3702--3716.e30},
	title = {Budding epithelial morphogenesis driven by cell-matrix versus cell-cell adhesion},
	volume = {184},
	year = {2021}}

@article{benguria1996variational,
	author = {Benguria, R. D. and Depassier, M. C.},
	journal = {Communications in Mathematical Physics},
	pages = {221--227},
	publisher = {Springer},
	title = {Variational characterization of the speed of propagation of fronts for the nonlinear diffusion equation},
	volume = {175},
	year = {1996}}

@article{benguria2004minimal,
	author = {Benguria, RD and Depassier, MC and M{\'e}ndez, V},
	journal = {Physical Review E},
	number = {3},
	pages = {031106},
	publisher = {APS},
	title = {Minimal speed of fronts of reaction--convection--diffusion equations},
	volume = {69},
	year = {2004}}

@article{benguria1994validity,
	author = {Benguria, R. D. and Depassier, M. C.},
	journal = {Physical Review Letters},
	number = {16},
	pages = {2272},
	publisher = {APS},
	title = {Validity of the linear speed selection mechanism for fronts of the nonlinear diffusion equation},
	volume = {73},
	year = {1994}}

@article{benguria1996speed,
	author = {Benguria, R. D. and Depassier, M. C.},
	journal = {Physical Review Letters},
	number = {6},
	pages = {1171},
	publisher = {APS},
	title = {Speed of fronts of the reaction--diffusion equation},
	volume = {77},
	year = {1996}}

@article{bruna2017diffusion,
	author = {Bruna, Maria and Chapman, S Jonathan and Robinson, Martin},
	journal = {SIAM Journal on Applied Mathematics},
	number = {6},
	pages = {2294--2316},
	publisher = {SIAM},
	title = {Diffusion of particles with short-range interactions},
	volume = {77},
	year = {2017}}

@article{bakerAspectRatio2,
	author = {Simpson, Matthew J. and Baker, Ruth E. and McCue, Scott W.},
	doi = {10.1103/PhysRevE.83.021901},
	journal = {Physical Review E},
	number = {2},
	numpages = {14},
	pages = {021901},
	publisher = {American Physical Society},
	title = {Models of collective cell spreading with variable cell aspect ratio: a motivation for degenerate diffusion models},
	volume = {83},
	year = {2011}}

@article{falco2022random,
	author = {Falc{\'o}, Carles},
	journal = {Physical Review E},
	number = {5},
	pages = {054103},
	publisher = {APS},
	title = {From random walks on networks to nonlinear diffusion},
	volume = {106},
	year = {2022}}

@article{Simpson2022ParameterIdentifiability,
	author = {Simpson, Matthew J. and Browning, Alexander P. and Warne, David J. and Maclaren, Oliver J. and Baker, Ruth E.},
	doi = {10.1016/j.jtbi.2021.110998},
	journal = {Journal of Theoretical Biology},
	pages = {110998},
	title = {Parameter identifiability and model selection for sigmoid population growth models},
	volume = {535},
	year = {2022}}

@article{Liu2024ParameterIdentifiabilityPDE,
	author = {Liu, Yue and Suh, Kevin and Maini, Philip K. and Cohen, Daniel J. and Baker, Ruth E.},
	doi = {10.1098/rsif.2023.0607},
	journal = {Journal of the Royal Society Interface},
	number = {212},
	pages = {20230607},
	title = {Parameter identifiability and model selection for partial differential equation models of cell invasion},
	volume = {21},
	year = {2024}}

@article{Alarcon2005MultipleScale,
	author = {Alarc{\'o}n, Tom{\'a}s and Byrne, Helen M. and Maini, Philip K.},
	doi = {10.1137/040603760},
	journal = {Multiscale Modeling \& Simulation},
	number = {2},
	pages = {440--475},
	title = {A multiple scale model for tumor growth},
	volume = {3},
	year = {2005}}

\end{document}